\documentclass[final]{svjour3}            
\smartqed  
\usepackage{hyperref}
\hypersetup{
    colorlinks=true,
    linkcolor=blue,
    filecolor=magenta,
    urlcolor=cyan,
    citecolor=black,
}
\usepackage{mathtools}
\usepackage[mathscr]{euscript}
\usepackage{graphicx}
 \usepackage{mathptmx}      
\usepackage{chicago-bibstyle}  
\usepackage[misc]{ifsym}
\usepackage{bbding}
\usepackage{appendix}
 \journalname{Bulletin of Mathematical Biology}
\begin{document}
\title{Complex Far-Field Geometries Determine the Stability of Solid Tumor Growth with Chemotaxis}



\author{Min-Jhe Lu$^1$ \and
        Chun Liu$^1$  \\
        John Lowengrub$^2$   \and
        Shuwang Li$^1$
}

\authorrunning{Min-Jhe Lu et al.} 

\institute{
\begin{itemize}
    \renewcommand{\labelitemi}{\Letter}
    \item  John Lowengrub \\
            lowengrb@math.uci.edu \\
   \item  Shuwang Li \\
            sli@math.iit.edu \\
    \renewcommand{\labelitemi}{$^1$}
    \item Department of Applied Mathematics,
          Illinois Institute of Technology, Chicago, IL 60616, USA\\
    \renewcommand{\labelitemi}{$^2$}
    \item Departments of Mathematics and Biomedical Engineering, Center for Complex Biological Systems, Chao Family Comprehensive Cancer Center, University of California at Irvine, Irvine, CA 92617, USA\end{itemize}
}

\date{Received: 12 December 2019 / Accepted: 26 February 2020}

\maketitle

\begin{abstract}
In this paper, we develop a sharp interface tumor growth model to study the effect of the tumor microenvironment using a complex far-field geometry that mimics a heterogeneous distribution of vasculature. Together with different nutrient uptake rates inside and outside the tumor, this introduces variability in spatial diffusion gradients. Linear stability analysis suggests that the uptake rate in the tumor microenvironment, together with chemotaxis, may induce unstable growth, especially when the nutrient gradients are large. We investigate the fully nonlinear dynamics using a spectrally accurate boundary integral method. Our nonlinear simulations reveal that vascular heterogeneity plays an important role in the development of morphological instabilities that range from fingering and chain-like morphologies to compact, plate-like shapes in two-dimensions.



\end{abstract}

\section{Introduction}
\label{sec:1}

 Tumor morphologies are determined by the competition between cell proliferation, cell death and chemotaxis, which are determined at least in part by the distribution of nutrients within the tumor and its microenvironment (e.g., \cite{cristini2005morphologic}). Thus, heterogeneities in vascular distributions and variable uptake rates of nutrients can lead to heterogeneous nutrient distributions that could induce diffusional instability through  nonuniform rates of cell proliferation, death and migration. Morphological instability, in turn,  is capable of bringing more available nutrients to the tumor by increasing its surface area to volume ratio. In particular, regions where instabilities first occur tend to grow at a faster rate than the rest of the tumor tissue  (\textit{e.g.}, differential growth) that further enhances the instabilities and leads to complicated tumor morphologies.

The process of nutrient diffusion and uptake inside tumors has been studied since the mid-1960s. See, for example, the reviews
(\cite{araujo2004history,fasano2006mathematical,roose2007mathematical,bellomo2008selected,lowengrub2009nonlinear,byrne2010dissecting,byrne2012mathematical,KimOthmer2015,alfonso2017biology,Yankeelov2018}) and books (\cite{cristini2010multiscale,cristini2017introduction}). Early biomechanical models of avascular tumor growth were proposed by Greenspan \cite{greenspan1976growth}. Later, Friedman and co-workers used analytical techniques to investigate the bifurcation behavior of the tumor growth and predict how the invasive morphology of a tumor develops \cite{friedman2001existence,friedman2001symmetry,friedman2006bifurcation,Hu20071,Hu20072,friedman2008stability}. Cristini et al. \cite{cristini2003nonlinear} pioneered the use of numerical simulations within the framework of a Darcy-law tumor growth model where the tumor/host interface is considered to be mathematically sharp. Along this direction, Pham et al.  \cite{pham2018nonlinear} developed a boundary integral method to simulate tumor growth using a Stokes flow model focusing on the viscosity ratio between the tumor and host tissues. Chemotaxis was introduced to the tumor model in \cite{cristini2009}, where a phase-field mixture model (see also, e.g., \cite{fritz2019local}) was used to study tumor growth in a simple microenvironmental geometry. A sharp interface model was derived from the phase field model using matched asymptotic expansions. Escher and Matioc \cite{escher2013analysis} analyzed local well-posedness and stability of this sharp interface model. Garcke et. al \cite{garcke2016cahn} derived a Cahn--Hilliard--Darcy model (see also \cite{garcke2018multiphase} for a multi-component variant) with chemotaxis and active transport using fundamental thermodynamic principles, the asymptotic limit of which serves as an extension of the previous sharp interface model to incorporate active transport. In this paper, we extend this work by developing, analyzing and numerically simulating a sharp interface model for tumor growth with chemotaxis in a complex far-field geometry that mimics a heterogeneous distribution of vasculature.

The main goals of this paper are to analyze the linear stability of the tumor model with chemotaxis in a complex microenvironment, and numerically simulate the fully nonlinear dynamics, using a novel boundary integral method. The work presented in this paper is unique in the following aspects. First, we consider a two-phase nutrient model where the uptake rates of the tumor and host tissues may be different, which significantly influences nutrient gradients that can drive morphological instabilities. Second, from the numerical perspective, we develop a new boundary integral method (BIM), which naturally incorporates the complex far-field geometries without approximation errors introduced by spatial meshes in the tumor/host domains (e.g., \cite{li2009solving}). The three coupled integral equations uniquely determine the nutrient, its normal gradient at the interface and the nutrient flux across the complex far-field boundary.  We believe that this is the first sharp interface model that solves coupled integral equations for the detailed nutrient distributions at the tumor/host interface and across the interface with spectral accuracy. Third, we investigate the significance of the microenvironmental variables (such as the vascular heterogeneity) on the fully nonlinear tumor dynamics.

The paper is organized as follows. In Sect. \ref{sec:2}, we formulate the two-phase nutrient field sharp interface model and non-dimensionalize the resulting systems. In Sect. \ref{sec:3}, we develop the BIM formulation. In Sect. \ref{sec:4}, we analyze the linear stability of the system. In Sect. \ref{sec:5}, we present our simulation results of the fully nonlinear system including a numerical convergence study, a comparison with linear analysis results and parameter studies under different far-field geometries. The conclusion is presented in Sect. \ref{sec:6}. In "Appendices A, B,C and D" we give a complete derivation of our linear stability analysis, the details of our numerical method including layer potential evaluations for boundary integrals, the small-scale decomposition to remove stiffness from the high-order derivatives in high curvature region on the interface and the semi-implicit time-stepping scheme to evolve the interface. The additional results on the effect of chemotaxis are shown in "Appendix E."
\section{Mathematical model}
\label{sec:2}

\begin{figure}
  \centering
  \includegraphics[width=0.5\textwidth]{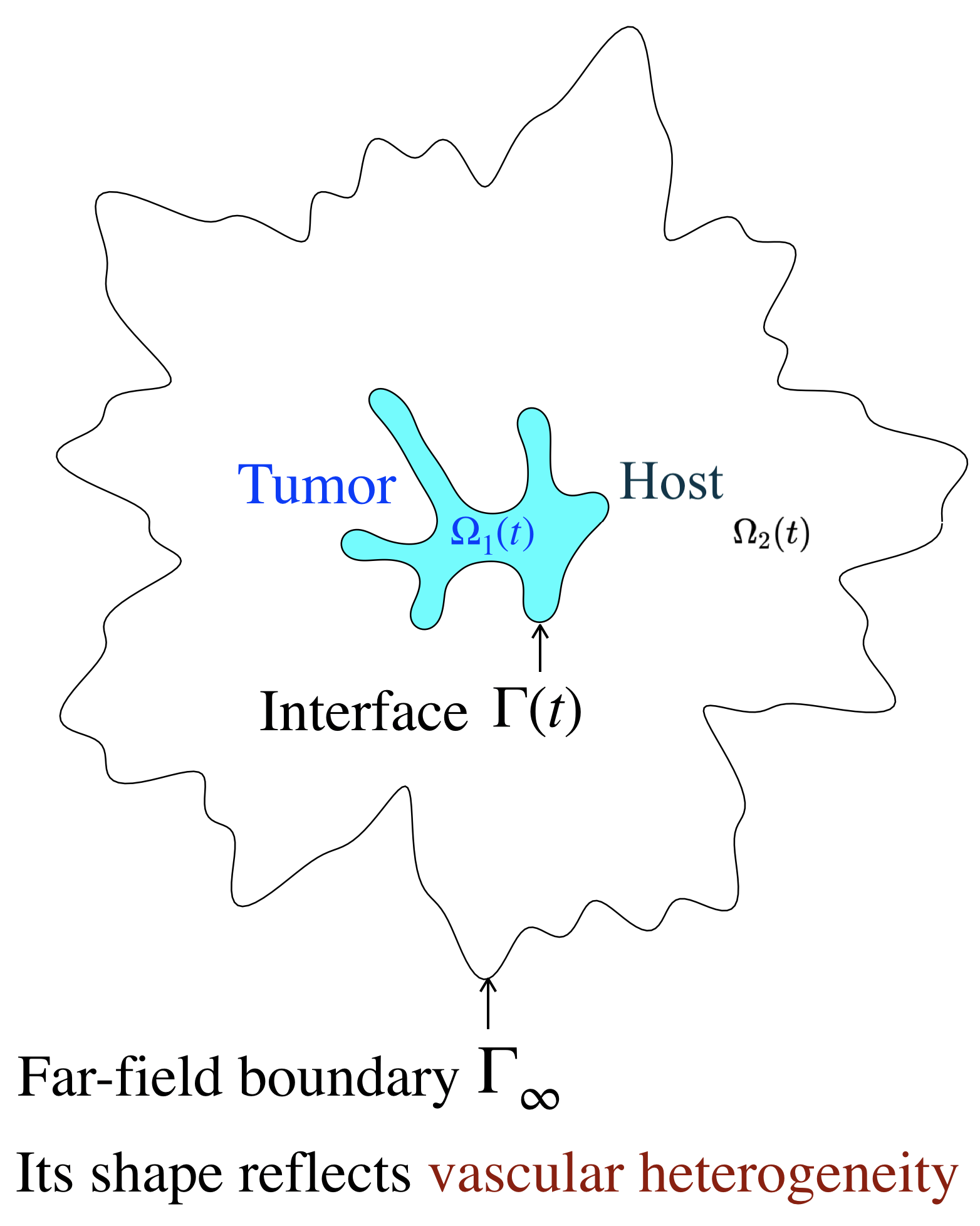}
\caption{Illustration of the computation domain of the tumor growth model with a complex far-field geometry.}
\label{fig:1}
\end{figure}

We follow the framework developed in  \cite{cristini2003nonlinear,cristini2009}. In particular, as illustrated in Fig. \ref{fig:1}, let $\Omega_1$ be the tumor tissue, $\Omega_2$ be the host tissue, $\Gamma(t)$ be the evolving tumor-host interface and $\Gamma_\infty$ be the fixed, far-field boundary the shape of which mimics the vascular heterogeneity. We assume that the soft tissue mechanics can be described using Darcy's law and that cell motion is driven by displacement due to nutrient-driven cell proliferation and death and chemotaxis of cells up nutrient gradients.

\subsection{A Sharp Interface Model with a Two-Phase Nutrient Field}
\label{sec:2.2}
\paragraph{Nutrient transport} Let $\sigma$  be the concentration of nutrients and growth-promoting factors; for simplicity we do not distinguish among the various biochemical molecules.
Then, the nutrient field in $\Omega_i$ satisfies the reaction-diffusion equation:
\begin{equation}\label{eq1}
\sigma_t=\nabla\cdot(D_i \nabla\sigma) - \lambda_i \sigma,
\end{equation}where $D_i,\lambda_i$ are the diffusion constants and uptake rates in $\Omega_i$ respectively. In general, \(D_i\) and $\lambda_i$ could be spatially and temporally heterogeneous or functions
 of $\sigma$ and $|\nabla\sigma|$. Because of the nutrient on the far-field boundary $\sigma^\infty$,
 this significantly limits the range of $\sigma$ in the domain, e.g., $0\le\sigma\le \sigma^\infty$, which should limit the
 effect of potential nonlinearities in $\sigma$. Although
 $|\nabla\sigma|$ can be large, we consider the simplest linear diffusion in which \(D_i\) and $\lambda_i$ are constants.
 The effects of nonlinearity and spatiotemporal heterogeneity can be considered in a future work.
Next, assuming the continuity of concentration and fluxes on $\Gamma(t)$:
\begin{align}
\left[\sigma\right]&\equiv\sigma_2-\sigma_1=0,\label{eq2}\\
\left[D_i\frac{\partial \sigma}{\partial \mathbf{n}}\right]&\equiv D_2\frac{\partial \sigma_2}{\partial \mathbf{n}}-D_1\frac{\partial \sigma_1}{\partial \mathbf{n}}=0,\label{eq3}
\end{align}
where $\sigma_i$ denotes $\displaystyle \lim_{\mathbf{x}_i\rightarrow\Gamma(t)}\sigma(\mathbf{x}_i)$ for $\mathbf{x}_i$ in $\Omega_i$ with $i=1,2$ and analogously for $\frac{\partial \sigma_i}{\partial\mathbf{n}}$.
The far-field boundary condition is:
\begin{equation}\label{eq4}
\sigma=\sigma^{\infty},~~~~{\rm on}~~\Gamma_\infty.
\end{equation}
\paragraph{Chemotaxis}
To model the directed migration of tumor cells up gradients of nutrients, we use the Chemo-Darcy law:
\begin{equation}\label{eq5}
\mathbf{u}=-\mu \nabla p + \chi_{\sigma}\nabla \sigma,~~~~{\rm in}~~\Omega_1,
\end{equation}
where $p$ is the pressure, $\mu$ is the cell mobility and $\chi_\sigma$ is the chemotaxis coefficient.
Then, assuming the tumor cell density is constant, the conservation of mass in the tumor is:
\begin{equation}
\nabla  \cdot \mathbf{u}
=
\lambda_M\frac{\sigma}{\sigma_\infty}-\lambda_A,~~{\rm in}~~  \Omega_1,
\end{equation}
where $\lambda_M,\lambda_A$ is the rate of mitosis and apoptosis, respectively. We assume there is no cell proliferation or death in the host tissue region, and for simplicity, we do not model any motion of the host tissue.
\paragraph{Boundary conditions}
The Laplace-Young condition for the internal pressure is assumed on $\Gamma$:
\begin{equation}\label{eq6}
(p)_{\Gamma}=\gamma \kappa,
\end{equation}
where $\gamma$ is a measure of cell-to-cell adhesion and $\kappa$ is the curvature of $\Gamma$. The
equation of motion for the interface $\Gamma$ is given by:
\begin{equation}\label{eq8}
V=-\mu \left(\frac{\partial p}{\partial \mathbf{n}}\right)_{\Gamma} + \chi_\sigma \left(\frac{\partial \sigma_1}{\partial \mathbf{n}}\right)_{\Gamma},
\end{equation}
where $\displaystyle\left(\frac{\partial \sigma_1}{\partial \mathbf{n}}\right)_{\Gamma}\equiv\lim_{\substack{\mathbf{x}_1\rightarrow\Gamma\\\mathbf{x}_1\in\Omega_1}}\left(\frac{\partial \sigma(\mathbf{x}_1)}{\partial \mathbf{n}}\right)_{\Gamma}$, and analogously for $p$.
  \subsection{Non-dimensionalization}
\label{sec:2.3}
We introduce the diffusion length $L$, the intrinsic taxis time scale $\lambda_\chi^{-1}$ and the characteristic pressure $p_s$ by:
\begin{equation}\label{eq9}
L=\sqrt{\frac{D_1}{\lambda_1}}, \quad
\lambda_\chi=\frac{\overline{\chi_{\sigma}}\sigma^{\infty}}{L^2}, \quad
p_s=\frac{\lambda_{\chi}L^2}{\mu},
\end{equation}
\noindent
where $\overline{\chi_\sigma}$ is a characteristic taxis coefficient.
The non-dimensional variables and parameters are:
\begin{equation}\label{eq10}
\widetilde{\sigma}=\frac{\sigma}{\sigma^{\infty}} , \quad
\widetilde{p}=\frac{p}{p_s}, \quad
\widetilde{\chi_{\sigma}}=\frac{\chi_{\sigma}}{\overline{\chi_{\sigma}}} , \quad
D=\frac{D_2}{D_1}, \quad
\lambda=\frac{\lambda_2}{\lambda_1}.
\end{equation}
\noindent
During the derivation, we have
\begin{align}
  {\lambda_\chi\widetilde\sigma_{\tilde t}}&=\lambda_1(\widetilde \Delta\widetilde\sigma-\widetilde\sigma),\\
  {\lambda_\chi\widetilde\sigma_{\tilde t}}&=\lambda_1(D\widetilde \Delta\widetilde\sigma-\lambda\widetilde\sigma),
\end{align}
Here, we assume the time scales of diffusion and uptake rate in the reaction diffusion process are similar, and since diffusion occurs faster than the taxis (e.g., minutes vs hours), we have $\lambda_{1} \gg \lambda_{\chi}$, which leads to quasi-steady reaction diffusion equations for nutrient field.
The dimensionless system is thus given by:
\paragraph{Modified Helmholtz Equations\\\\}

 In $\Omega_i$:
\begin{equation}\label{modifiedHH}
\widetilde{\Delta}\widetilde{\sigma}_i-\mu_i^2\widetilde{\sigma}_i=0,
\end{equation}where $\mu_1=1, \mu_2=\frac{1}{\Lambda}$ and $\Lambda=\sqrt{{\frac{D}{\lambda}}}$ is the relative diffusional penetration length.
\paragraph{Non-dimensional Chemo-Darcy's law\\\\}
 In $\Omega_1$:
\begin{equation}\label{eq14}
\mathbf{\widetilde{u}}=-\widetilde{\nabla} \widetilde{p} + \widetilde{\chi_{\sigma}} \widetilde{\nabla} \widetilde{\sigma}.
\end{equation}

\paragraph{Conservation of tumor mass\\\\}
 In $\Omega_1$:
\begin{equation}\label{eq16}
\widetilde{\nabla} \cdot \widetilde{\mathbf{u}}=\mathcal{P} \widetilde{\sigma}- \mathcal{A},
\end{equation}where $\mathcal{P}=\frac{\lambda_M}{\lambda_\chi}$ represents the relative rate of cell mitosis to taxis and $\mathcal{A}=\frac{\lambda_A}{\lambda_\chi}$ represents the relative rate of cell apoptosis to taxis.
\paragraph{Boundary Conditions\\\\}

On $\Gamma$:
\begin{equation}\label{eq18}
\widetilde{\sigma}_1 =\widetilde{\sigma}_2.
\end{equation}\begin{equation}\label{eq19}
\frac{\widetilde{\partial} \widetilde{\sigma}_1}{\widetilde{\partial} \widetilde{\mathbf{n}}}=D\frac{\widetilde{\partial} \widetilde{\sigma}_2}{\widetilde{\partial} \widetilde{\mathbf{n}}}.
\end{equation}
\begin{equation}\label{eq20}
\widetilde{p}=\widetilde{\mathcal{G}}^{-1} \widetilde{\kappa},
\end{equation}where $\widetilde{\mathcal{G}}^{-1}
=\frac{\mu \gamma}{\lambda_\chi L^3}\text{ represents the relative strength of cell--cell interactions (adhesion)}.$\\\\
On $\Gamma_\infty$:
\begin{equation}\label{eq22}
\widetilde{\sigma}=1.
\end{equation}
\paragraph{Equation of motion of $\Gamma$\\\\}
On $\Gamma$:
\begin{equation}
\widetilde{V}=-\frac{\widetilde{\partial}\widetilde{p}}{\widetilde{\partial}\widetilde{\mathbf{n}}}
+\widetilde{\chi_{\sigma}}\frac{\widetilde{\partial}\widetilde{\sigma}}{\widetilde{\partial}\widetilde{\mathbf{n}}}.
\end{equation}

\section{Boundary Integral Method (BIM) Reformulation}\label{sec:3}

From Eqs.(\ref{modifiedHH}), (\ref{eq14}) and (\ref{eq16}), we have the Poisson equation for the non-dimensional pressure $\widetilde{p}$:
\begin{equation}\label{eq23}
-\widetilde{\Delta} \widetilde{p} = (\mathcal{P}-\widetilde{\chi_{\sigma}})\widetilde{\sigma}-\mathcal{A}.
\end{equation}
Using the algebraic transformation from \cite{cristini2003nonlinear} to define the modified pressure:
\begin{equation}\label{eq24}
\overline{p}=\widetilde{p}+(\mathcal{P}-\widetilde{\chi_{\sigma}})\widetilde{\sigma} - \mathcal{A}\frac{\mathbf{\widetilde{x}}\cdot \mathbf{\widetilde{x}}}{2 d},
\end{equation}
where $d$ is the dimension of $\mathbf{R}^d \supseteq \Omega_i$.
Then, from Eq. (\ref{modifiedHH}) we obtain
\begin{equation}\label{eq25}
-\widetilde{\Delta} \overline{p}= -\widetilde{\Delta} \widetilde{p} - (\mathcal{P}-\widetilde{\chi_{\sigma}})\widetilde{\sigma}+\mathcal{A}=0.
\end{equation}\\
Dropping all tildes and overbars for brevity, the non-dimensional system of equations is as follows. We have the Laplace equation for the modified pressure $p$:\\\\
In $\Omega_1$:
\begin{equation}\label{eq26}
\Delta p=0.
\end{equation}
On $\Gamma$:
\begin{equation}\label{eq27}
p=\mathcal{G}^{-1} \kappa
+(\mathcal{P}-\chi_{\sigma})\sigma - \mathcal{A}\frac{\mathbf{x}\cdot \mathbf{x}}{2 d},
\end{equation}
where $\sigma$ satisfies the two-phase modified Helmholtz equations:\\\\
In $\Omega_i$:
\begin{equation}\label{eq28}
\Delta \sigma_i - \mu_i^2\sigma_i =0.
\end{equation}
On $\Gamma$:
\begin{equation}\label{eq29}
\sigma_1=\sigma_2,
\end{equation}
\begin{equation}\label{eq30}
\frac{\partial \sigma_1}{\partial \mathbf{n}}
=D\frac{\partial \sigma_2}{\partial \mathbf{n}}.
\end{equation}
On $\Gamma_\infty$:
\begin{equation}\label{eq31}
\sigma_2=1.
\end{equation}
The equation of motion for the interface is then:
\begin{equation}\label{eq32}
V=-\frac{\partial p}{\partial \mathbf{n}}
+ \mathcal{P} \frac{\partial \sigma}{\partial \mathbf{n}}
-\mathcal{A} \frac{\mathbf{n}\cdot \mathbf{x}}{d}.
\end{equation}From potential theory, the solutions of Eqs. \eqref{eq26}, \eqref{eq28} can be formulated as boundary integrals with single-layer and double-layer potentials. We use a direct formulation for $\sigma$ and an indirect formulation for $p$.

\subsection{Direct BIM Formulation for the Two-Phase Modified Helmholtz Equations}
Consider the Green's function for modified Helmholtz equations in $\Omega_i$:
\begin{equation}\label{eq33}
\Delta G_i -\mu_i^2 G_i = -\delta_{\mathbf{x_i}},
\end{equation}where $G_i=G_i(\mathbf{x_i},\mathbf{x}')$, $\mathbf{x_i}\in \Omega_i$ is the source point, $\mathbf{x}'$ is the field point, and $\delta_{\mathbf{x_i}}(\mathbf{x_i},\mathbf{x}')$ is the Dirac delta function. The fundamental solution of Eq. \eqref{eq33} is
\begin{equation}\label{eq34}
G_i(\mathbf{x_i},\mathbf{x}')=\frac{1}{2 \pi} K_0(\mu_i r),
\end{equation}
where $K_0$ is a modified Bessel function of the second kind and $r=|\mathbf{x_i}-\mathbf{x}'|$.
Multiplying  Eq. \eqref{eq28} by $G_i$ , Eq. \eqref{eq33} by $-\sigma_i$ and summing them up, we get
\begin{equation}\label{eq35}
\sigma_i \delta_{\mathbf{x_i}}=G_i \Delta \sigma_i - \sigma_i \Delta G_i.
\end{equation}
Integrating Eq. \eqref{eq35} over $\Omega_i$ and using Green's second identity, we have
\begin{equation}\label{eq36}
\sigma_i(\mathbf{x_i}) =\int_{\Gamma_i}
\left(G_i \frac{\partial\sigma_i}{\partial\mathbf{n_i}'}
-\sigma_i \frac{\partial G_i}{\partial \mathbf{n_i}'}\right) ds',
\end{equation}where $\Gamma_1=\Gamma$, $\Gamma_2=\Gamma \cup \Gamma_{\infty}$, $\mathbf{n}_1'=\mathbf{n}' $ is the unit outer normal on $\Gamma$, $\mathbf{n_2}'=-\mathbf{n}'$ on $\Gamma$ and $\mathbf{n_2}'= \mathbf{n}'$ on $\Gamma_{\infty}$.\\\\
Following \cite{kress2013linear}, letting $\mathbf{x_1} \to \mathbf{x} \in \Gamma, \mathbf{x_2} \to \mathbf{x} \in \Gamma$ and $\mathbf{x_2} \to \mathbf{x_{\infty}} \in \Gamma_{\infty}$ in Eq. \eqref{eq36} and using Eqs. \eqref{eq29}, \eqref{eq30}, \eqref{eq31}, we obtain the following coupled set of boundary integral equations:

\paragraph{On $\Gamma$:}
\begin{equation}\label{eq40}
\frac{1}{2} \sigma_1 +
\int_{\Gamma}
\left(
\sigma_1
\frac{\partial G_1}{\partial \mathbf{n}'}
- G_1
\frac{\partial \sigma_1}{\partial \mathbf{n}'}
\right) ds' = 0,
\end{equation}\\
\begin{equation}\label{eq41}
-\frac{1}{2} \sigma_1+
\int_{\Gamma}
\left(
\sigma_1
\frac{\partial G_2 }{\partial \mathbf{n}'}
- G_2
\frac{1}{D} \frac{\partial \sigma_1}{\partial \mathbf{n}'}
\right) ds'
+\int_{\Gamma_\infty}
G_2
 \frac{\partial \sigma_2 }{\partial \mathbf{n_\infty}'}
 ds'
=
\int_{\Gamma_\infty}
\frac{\partial G_2}{\partial \mathbf{n_\infty}'}
ds'.
\end{equation}
\paragraph{On $\Gamma_\infty$:}
\begin{equation}\label{eq42}
\int_{\Gamma}
\left(
\sigma_1
\frac{\partial G_2 }{\partial \mathbf{n}'}
- G_2
\frac{1}{D} \frac{\partial \sigma_1}{\partial \mathbf{n}'}
\right) ds'
+\int_{\Gamma_\infty}
G_2
 \frac{\partial \sigma_2}{\partial \mathbf{n_\infty}'}
ds'= \frac{1}{2}+\int_{\Gamma_\infty} \frac{\partial G_2 }{\partial \mathbf{n_\infty}'}ds'.
\end{equation}This system is solved for the three unknowns $\sigma_1, \frac{\partial \sigma_1}{\partial \mathbf{n}'}$ on $\Gamma(t)$ and $\frac{\partial \sigma_2}{\partial \mathbf{n}'_\infty}$ on $\Gamma_\infty$. For details on how these integrals are discretized and how the system is solved, we refer the reader to "Appendix B."

\subsection{Indirect BIM Formulation for Laplace Equation}
Consider the Green's function for the Laplace equation in $\Omega_1$:
\begin{equation}\label{eq43}
\Delta G  = \delta_{\mathbf{x}},
\end{equation}where $G=G(\mathbf{x},\mathbf{x}')$, $\mathbf{x}\in \Omega_1$ is the source point, $\mathbf{x}'$ is the field point and $\delta_{\mathbf{x}}(\mathbf{x},\mathbf{x}')$ is the Dirac functional. The fundamental solution for Eq. \eqref{eq43} is
\begin{equation}\label{eq44}
G(\mathbf{x},\mathbf{x}')=\frac{1}{2 \pi} \ln r,
\end{equation}where $r=|\mathbf{x}-\mathbf{x}'|$.
From potential theory \cite{kress2013linear}, the modified pressure $p$ can be formulated as the double-layer potential defined by:
\begin{equation}\label{eq46}
p(\mathbf{x})=
    (\mathscr{D}\eta) (\mathbf{x})\equiv
    \int_{\Gamma}
    \frac{\partial G(\mathbf{x},\mathbf{x}')}{\partial \mathbf{n}'}
    \eta(\mathbf{x}')ds',
\end{equation}
where $\mathbf{n}'$ is the unit outer normal on $\Gamma$. Using the properties of the double-layer potential, we may recast Eqs. \eqref{eq26}, \eqref{eq27} into a second-kind Fredholm integral equation with the unknown $\eta$ on $\Gamma$:
\begin{equation}\label{eq47}
    (-\frac{1}{2}+\mathscr{D})(\eta)=
    \mathcal{G}^{-1} \kappa
+(\mathcal{P}-\chi_{\sigma})\sigma - \mathcal{A}\frac{\mathbf{x}\cdot \mathbf{x}}{2 d}.
\end{equation}The term $\frac{\partial p}{\partial \mathbf{n}}$ (normal derivative of double-layer potential) needed for evaluating normal velocity V in Eq. \eqref{eq32} can be computed as \cite{kress2013linear}:
\begin{equation}\label{eq49}
    \frac{\partial p}{\partial \mathbf{n}}
    =\frac{d}{ds} \mathscr{S}(\eta'_{s'}),
\end{equation}
where $s$ is the arclength parameter and $\mathscr{S}(\eta)$ is the single layer potential defined as:
\begin{equation}\label{eq45}
    (\mathscr{S}\eta)(\mathbf{x})\equiv\int_{\Gamma}G(\mathbf{x},\mathbf{x}')\eta(\mathbf{x}')ds'.
\end{equation}
"Appendices C and D" describe the corresponding discretizations used and the time-stepping method for evolving the interface $\Gamma(t)$ in time.

\section{Linear Analysis}
\label{sec:4}
In this section we present the results of a linear stability analysis (details are provided in "Appendix A") of the non-dimensional sharp interface equations \eqref{eq26}--\eqref{eq32} reformulated in the preceding section. The linear stability of perturbed radially symmetric tumors was previously analyzed in \cite{cristini2003nonlinear,cristini2009,garcke2016cahn}. Here, we extend their results to take into account the two-phase nutrient field in two dimensions.
Consider a $l^{th}$ mode perturbation of a radially symmetric tumor interface $\Gamma$:
\begin{equation}
r(\theta,t)=R(t)+\delta(t) \cos l \theta,
\end{equation}
where $r$ is the tumor/host interface, $R$ is the radius of the underlying circle, $\delta$ is the dimensionless perturbation size and $\theta$ is
the polar angle. We deduce that the evolution equation for the tumor radius $R$ is given by:
\begin{equation}\label{Rdot}
\frac{d{R}}{dt}=C \mathcal{P}-\frac{\mathcal{A} R}{d},
\end{equation}
where $C=\mu_{1} A_{1} I_{1}\left(\mu_{1} R\right)$.
The equation of the shape perturbation $\frac{\delta}{R}$ is given by:
\begin{eqnarray}{}
\left(\frac{\delta}{R}\right)^{-1}
\frac{d}{dt}\left(\frac{\delta}{R}\right)
&=&-\mathcal{G}^{-1} \frac{l\left(l^{2}-1\right)}{R^{3}}+\frac{l}{d} \mathcal{A}\nonumber\\
&&+\mathcal{P}\left(\mu_{1}^{2} A_{1} I_{0}\left(\mu_{1} R\right)+B_{1}\left(\mu_{1} I_{l-1}\left(\mu_{1} R\right)-\frac{l}{R} I_{l}\left(\mu_{1} R\right)\right)-\frac{2}{R} {C}\right)\nonumber\\
&&-\left(\mathcal{P}-\chi_{\sigma}\right)\left(B_{1} \frac{l}{R} I_{l}\left(\mu_{1} R\right)+\frac{l}{R} C\right),
\end{eqnarray}
where $I_{l-1}(R)$ and $I_{l}(R)$ are modified Bessel functions of the first kind. A complete derivation and the expressions of $A_1, B_1$, which contain the parameters $D$, $R_{\infty}$ and $\lambda$ are given in "Appendix A." We remark that when $\lambda=0$, we recover the linear stability analysis in \cite{cristini2009}. Furthermore, when the chemotaxis is absent and the diffusion constant outside the tumor is large, we recover the linear stability analysis in \cite{cristini2003nonlinear}.

We next characterize the stability regime by determining $\mathcal{A}=\mathcal{A}_c$ as a function of the unperturbed radius $R$ such that $\frac{d}{dt}{\left(\frac{\delta}{R}\right)}=0$:
\begin{eqnarray}{} \label{apoptc}
\mathcal{A}_{c}= \mathcal{G}^{-1} \frac{d\left(l^{2}-1\right)}{R^{3}}
&&-\mathcal{P} \frac{d}{l}\left(\mu_{1}^{2} A_{1} I_{0}\left(\mu_{1} R\right)+B_{1}\left(\mu_{1} I_{l-1}\left(\mu_{1} R\right)-\frac{l}{R} I_{l}\left(\mu_{1} R\right)\right)-\frac{2}{R} {C}\right)\nonumber\\ &&+\left(\mathcal{P}-\chi_{\sigma}\right) \frac{d}{l}\left(B_{1} \frac{l}{R} I_{l}\left(\mu_{1} R\right)+\frac{l}{R} {C}\right).
\end{eqnarray}
Here, $\mathcal{A}_c$ is the critical value of apoptosis that divides regimes of stable growth ($\mathcal{A} < \mathcal{A}_c$, \textit{e.g.}, below the curve) and regimes of unstable growth ($\mathcal{A} > \mathcal{A}_c$, \textit{e.g.}, above the curve) for a given mode $l$ and penetration length $\Lambda$. We focus on the parameters $\mathcal{G}^{-1}=0.001,\mathcal{P}=0.5$.

\begin{figure}
\centering
\includegraphics[width=\textwidth]{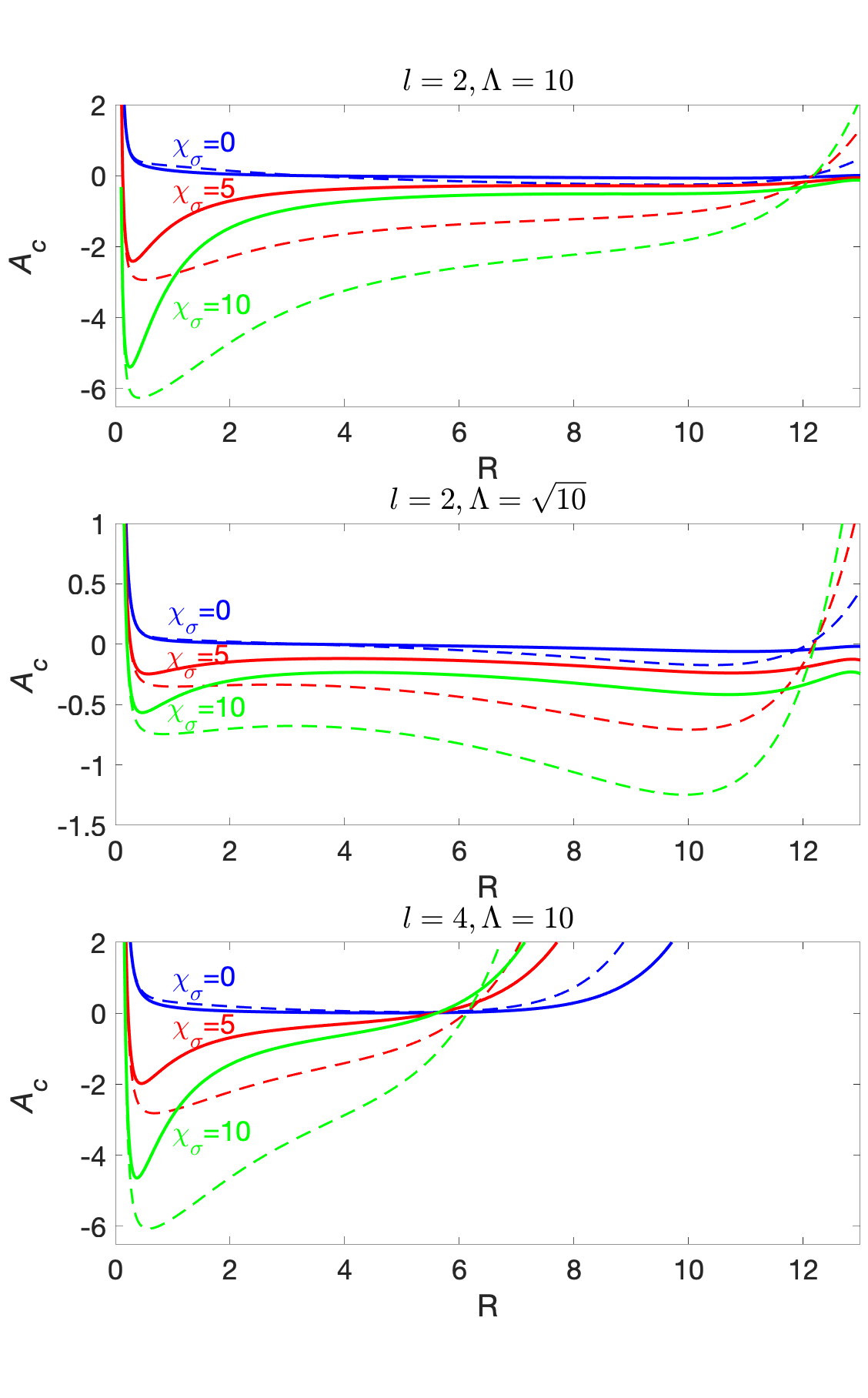}
\caption{Critical apoptosis parameter $\mathcal{A}_c$ as a function of unperturbed radius $R$ from equation \eqref{apoptc}, far-field-boundary $R_\infty= 13$, and $\chi_\sigma$ labeled. Solid: $D = 1$; Dashed: $D = 100$. See text for details.}\label{fig:apoptc}
\end{figure}

In Fig. \ref{fig:apoptc}, we plot $\mathcal{A}_c$ as a function of $R$, $D=1$ (solid), $D=100$ (dashed), and $\chi_\sigma$ as labeled. The top figure with $l=2,\Lambda=10$ reveals that the unstable regime expands with stronger taxis and most of the solid curves are pulled downward under richer supply of nutrient $(D=100)$.
Notice that we can also see the taxis stabilizes the dynamics when the tumor grows close to far-field boundary (here $R_\infty=13$) because the gradients are decreased. The middle figure with $l=2,\Lambda=\sqrt{10}$ shows a similar profile although the curves shift upward compared to those when $\Lambda=10$. Thus, the smaller the penetration length (or equivalently the larger the uptake rate in the host tissue), the more stable the dynamics becomes. In the bottom figure with $l=4, \Lambda=10$ we see that the unstable regime is narrowed, which indicates that mode $l = 4$ is more stable for this set of parameters.

\section{Results}
\label{sec:5}
Now we investigate the nonlinear evolution of tumor proliferation and invasion into surrounding host tissue. The details of the numerical methods are provided in "Appendices B-D". "Appendix B" provides details on the discretization and solution of the boundary integral equations (\ref{eq40})-(\ref{eq42}), and the methods for evolving the tumor/host interface $\Gamma$ in time are given in "Appendices C and D." In addition to the effects of $D$ and $\chi_\sigma$, investigated previously in  \cite{macklin2007nonlinear,cristini2009,garcke2016cahn}, we also consider the influence of the uptake-ratio $\lambda$ and complex-shaped, far-field geometries $\Gamma_\infty$ to mimic the effects of vascular heterogeneities. We fix the relative diffusional penetration length $\Lambda=10$, the non-dimensional proliferation and adhesion parameters $\mathcal{P}=0.5$, $\mathcal{G}^{-1}=0.001$, unless otherwise specified.

The final time snapshots presented from Figures $\ref{fig:5}$ to $\ref{fig:6-2}$ are chosen close to the time step which our numerical simulation stops and do not correspond to steady states. The computation terminates due to a lack of numerical resolution because the curvature of the tumor boundary becomes very large. To proceed significantly farther in time, many more grid points along the interface are needed to resolve the nearly singular curvature.

\subsection{Numerical Convergence in Time and Space}
\label{sec:4.1}
\begin{figure}
\begin{minipage}{1.0\linewidth}
\centering
\includegraphics[width=\textwidth]{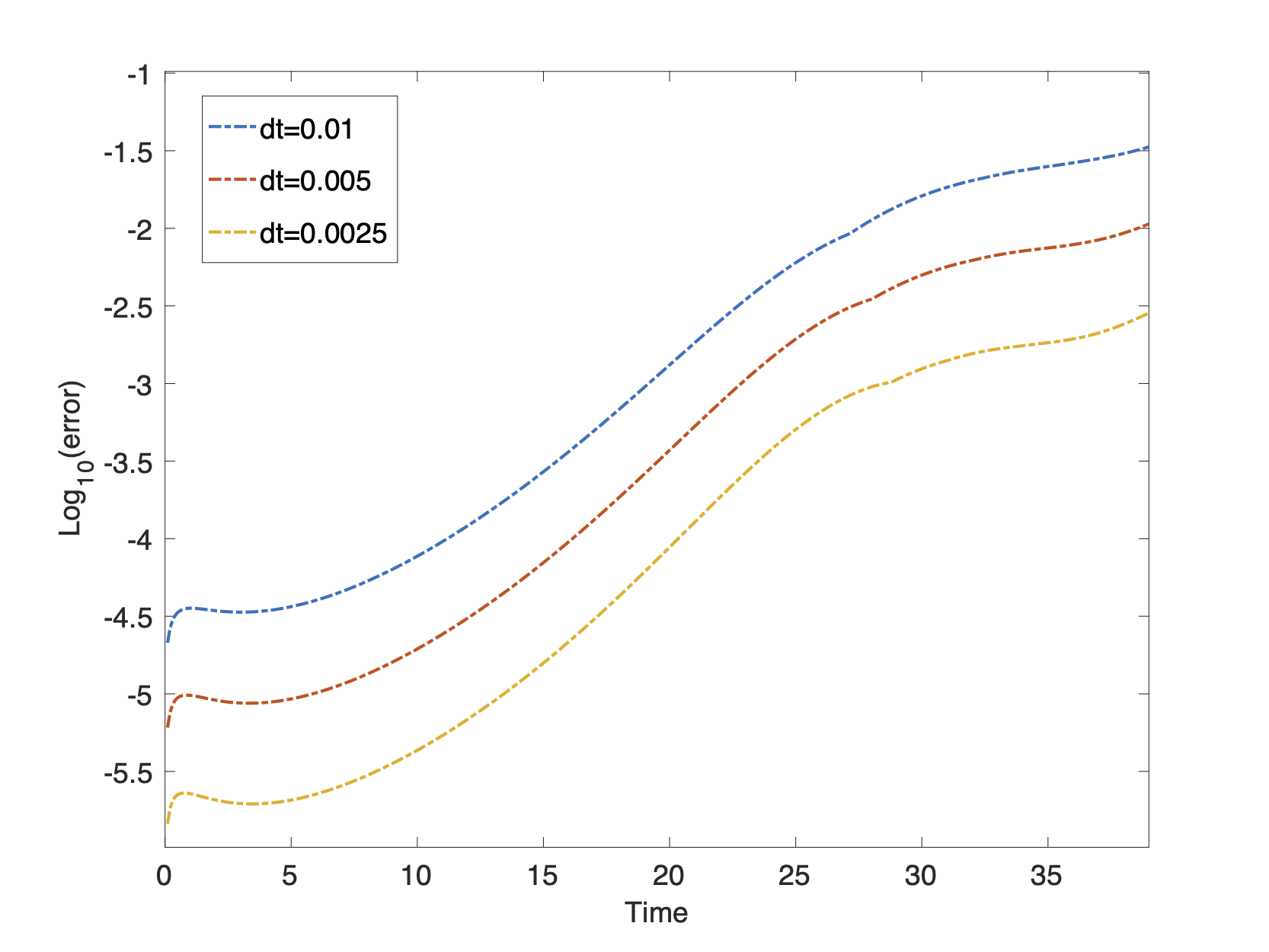}[a]
\end{minipage}
\begin{minipage}{1.0\linewidth}
\centering
\includegraphics[width=\textwidth]{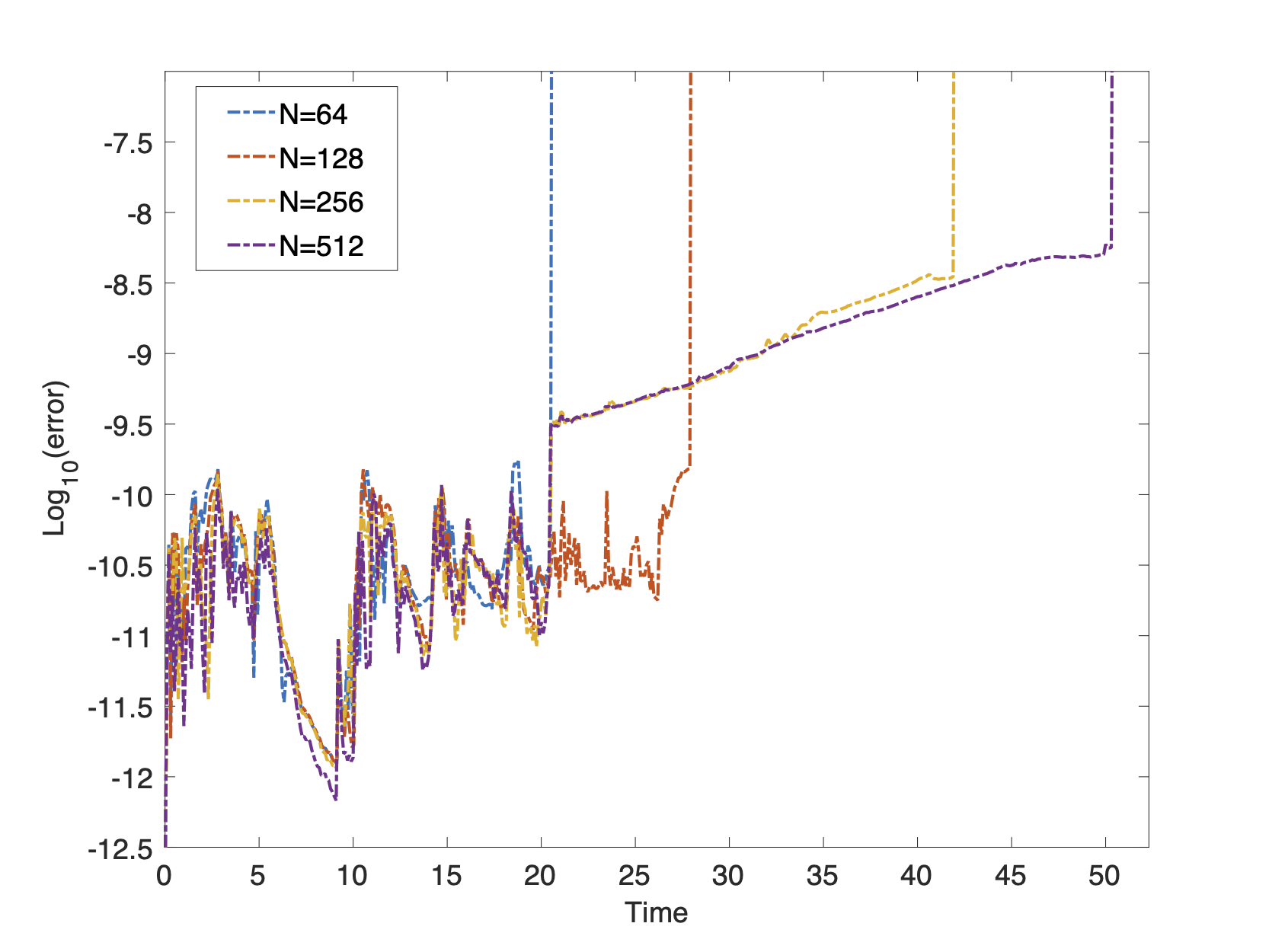}
\end{minipage}
\begin{minipage}{1.0\linewidth}
\centering
\includegraphics[width=\textwidth]{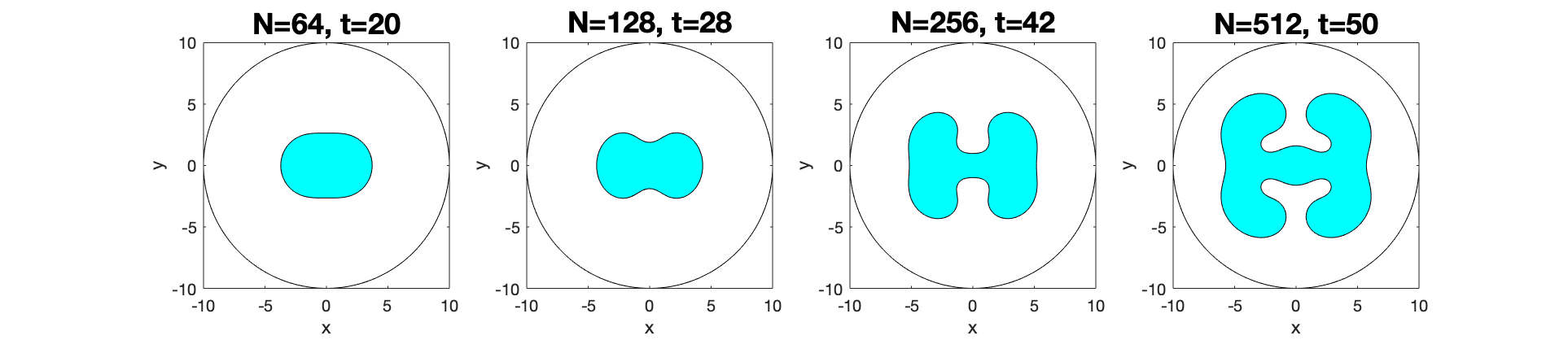}[b]
\end{minipage}

\caption{Temporal [a] and spatial [b] resolution studies for the case with parameters $D=1000$, $\lambda=1$, $\chi_{\sigma}=10$, $\mathcal{P}=1$, $\mathcal{A}=0.5$,  $\mathcal{G}^{-1}=0.05$. The far-field boundary is circular with radius $R_\infty=10$, and initial tumor boundary is $r=2.0+ 0.1\cos (2 \theta)$. In [a], the errors shown are calculated as the maximal differences
between the solution with $\Delta t=10^{-3}$ and those with $\Delta t=2.5\times10^{-3},5\times10^{-3},1\times10^{-2}$ up to time $t=39$; In [b] the errors are calculated as the maximal differences
between the solution with $N=1024$ and those with $N=64, 128, 256, 512$ up to their last time step $t=20, 28, 42, 50$. The corresponding tumor morphologies are shown in the last row.
 }
\label{fig:2}
\end{figure}

In this section, we test the convergence of our method.
First, we present a temporal resolution study. The errors from the time discretization are shown in Fig. \ref{fig:2}[a]. The maximal differences on the tumor interface between the simulation with $\Delta t=1\times 10^{-3}$ and $\Delta t=2.5\times10^{-3},5\times10^{-3},1\times10^{-2}$ , respectively,  are plotted versus time. In all cases, the number of spatial collocation points is $N = 256$. An improvement in accuracy with a factor of 4 is observed when $\Delta t$ is halved, indicating a second order convergence rate. This is expected since the time stepping scheme is second order accurate ("Appendix D").

In space, the accuracy of our simulation is established by a resolution study of the simulation shown in the inset of Fig. \ref{fig:2}[b]. The spatial error is investigated by varying the number $N$ of spatial collocation points representing the tumor interface. The maximal differences on the tumor interface between the simulation with $N=1024$ and those with $N=64,128,256,512$ respectively are plotted versus time. In all cases the time step is $\Delta t=1\times 10^{-2}.$ At early times, the error is dominated by the tolerance for solving the integral equations $(1\times 10^{-10}).$ This is consistent with the spectral accuracy of our method. The calculations with smaller $N$ cease at much earlier times than the larger $N$ due to the failure in solving the integral equations with the given tolerance at the spatial resolution dictated by $N$. In particular, the solution computed with $N=512$ is accurate to within $1\times 10^{-8.5}$, until about $t\approx 50$, where a larger $N$  is needed to resolve the regions of the interface that are in near contact.

\subsection{Comparison with Linear Analysis When the Far-Field Geometry is Circular}
\label{sec:4.2}
In this section, we compare the nonlinear simulation with the linear stability analysis results from Sec. \ref{sec:4} using a circular far-field boundary.

  \begin{figure}
  \begin{minipage}{1.0\linewidth}
  \centering
     \includegraphics[width=\textwidth]{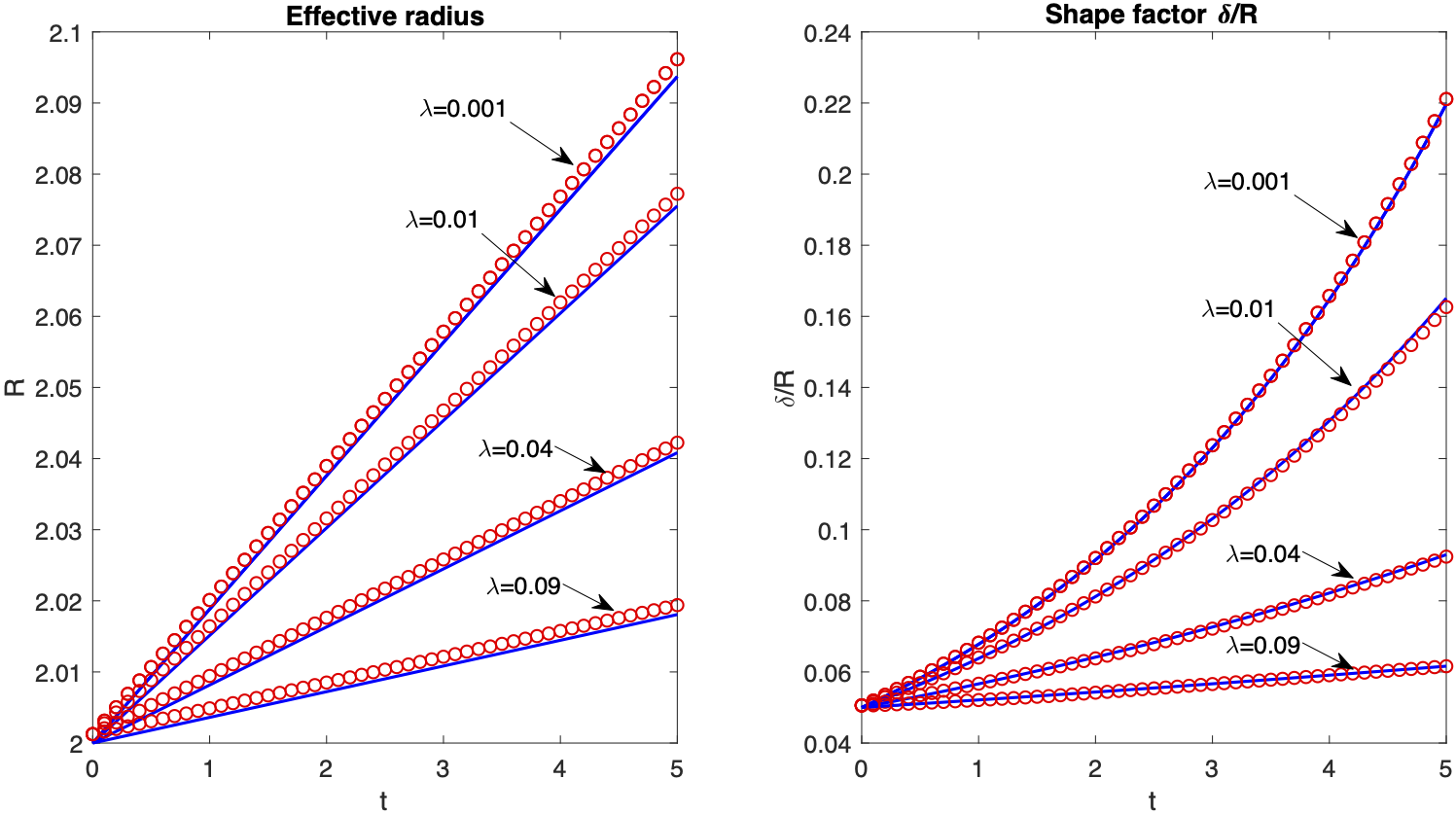}[a]
  \end{minipage}
\begin{minipage}{1.0\linewidth}
\centering
     \includegraphics[width=\textwidth]{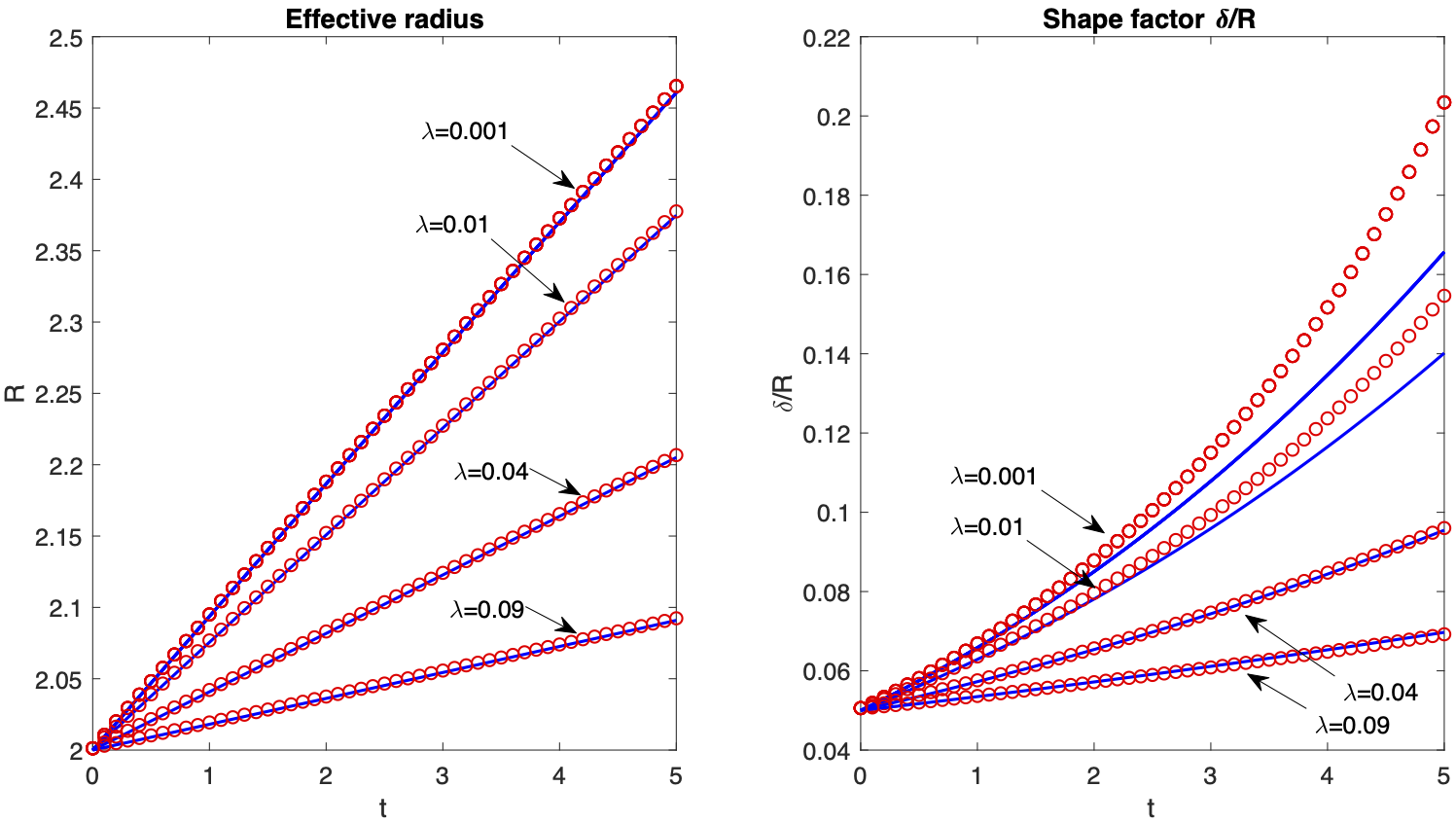}[b]
   	 \end{minipage}
\begin{minipage}{1.0\linewidth}
\centering
     \includegraphics[width=\textwidth]{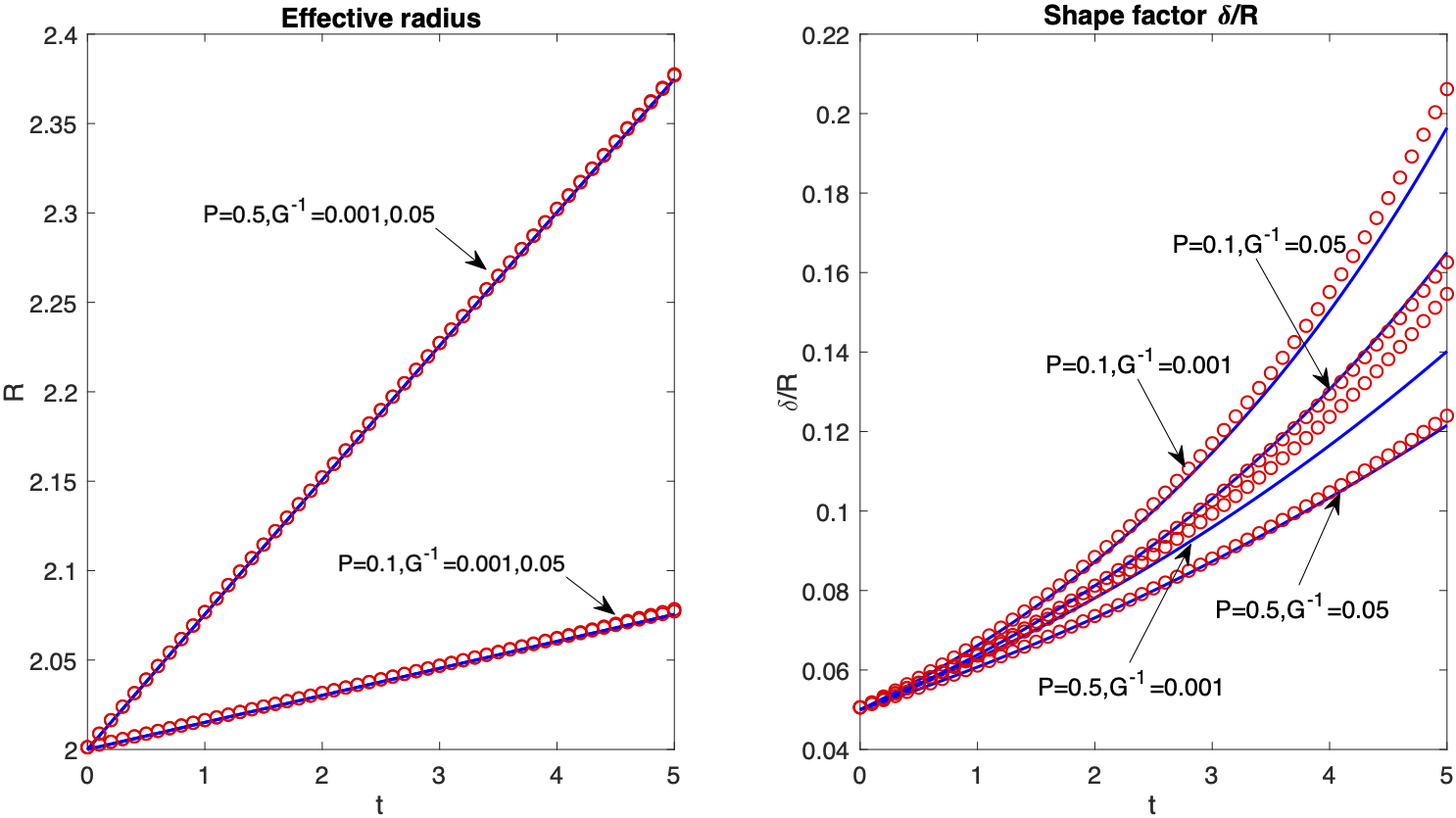}[c]
   	 \end{minipage}
     \caption{Time evolution of the effective radii $R$ and shape factors $\delta/R$. Blue lines: linear solutions; Red circles: nonlinear simulations. In [a] and [b]: $\lambda=0.001,0.01,0.04,0.09$; and in [c]: $\lambda=0.01$. The parameters are $D=1$, $\chi_\sigma=5$,  $\mathcal{P}= 0.1$ in [a] and $\mathcal{P}=0.5$ in [b], $\mathcal{A}=0$, $\mathcal{G}^{-1}=0.05$ in [a] and $\mathcal{G}^{-1}= 0.001$ in [b]. The far-field boundary is a circle with radius $R_\infty=13$, and the initial tumor surface is $r=2.0+0.1 \cos (2 \theta)$. Here, $N=64$ and $\Delta t=0.05.$ See text for details.}
     \label{fig:3}
  \end{figure}
In Fig. \ref{fig:3}[a] and [b], the linear stability analysis results (blue lines) for the effective radius $R$ (left) and shape factor (right) are compared with the nonlinear simulations (red circles) to show the stabilizing effect of the parameter $\lambda$ in the nutrient-poor regime. The nonlinear shape factor is calculated by ${\rm max}_\Gamma | \left(|\mathbf{x}|/R-1\right)|$.  The effective radius is defined as $R=\sqrt{\frac{\text{Area}}{\pi}}$. Note that the case with $\lambda=0.001\approx 0$ approximates the single-phase nutrient model used in \cite{cristini2009} and others.

In Fig. \ref{fig:3}[a] and [b], the larger the uptake-ratio $\lambda$ is, the slower and stabler the tumor grows. This is because the nutrient gradients decrease as $\lambda$ increases, which agrees with what we found in Fig. \ref{fig:apoptc}. Further in Fig. \ref{fig:3}[a], we see with $\mathcal{P}=0.1, \mathcal{G}^{-1}=0.05$, there is very good agreement between the linear solutions and nonlinear simulations. However, in Fig. \ref{fig:3}[b] when $\mathcal{P}=0.5, \mathcal{G}^{-1}=0.001$, while there is very good agreement between the linear and nonlinear effective radii, the linear solutions under predict the shape factors of the nonlinear results when $\lambda=0.001,\ 0.01$ at long times. In other words, nonlinearity enhances the morphological instability.

In Fig. \ref{fig:3}[c], first we notice on the left that the tumor with larger $\mathcal{P}$ grows at a faster rate, as indicated by Eq. \eqref{Rdot} and is independent of $\mathcal{G}^{-1}$. Second, we notice that tumors with smaller $\mathcal{P}$ and $\mathcal{G}^{-1}$ are more unstable and as before linear theory under predicts the instability.

\begin{figure}
\begin{minipage}{1.0\linewidth}
\centering
\includegraphics[width=\textwidth]{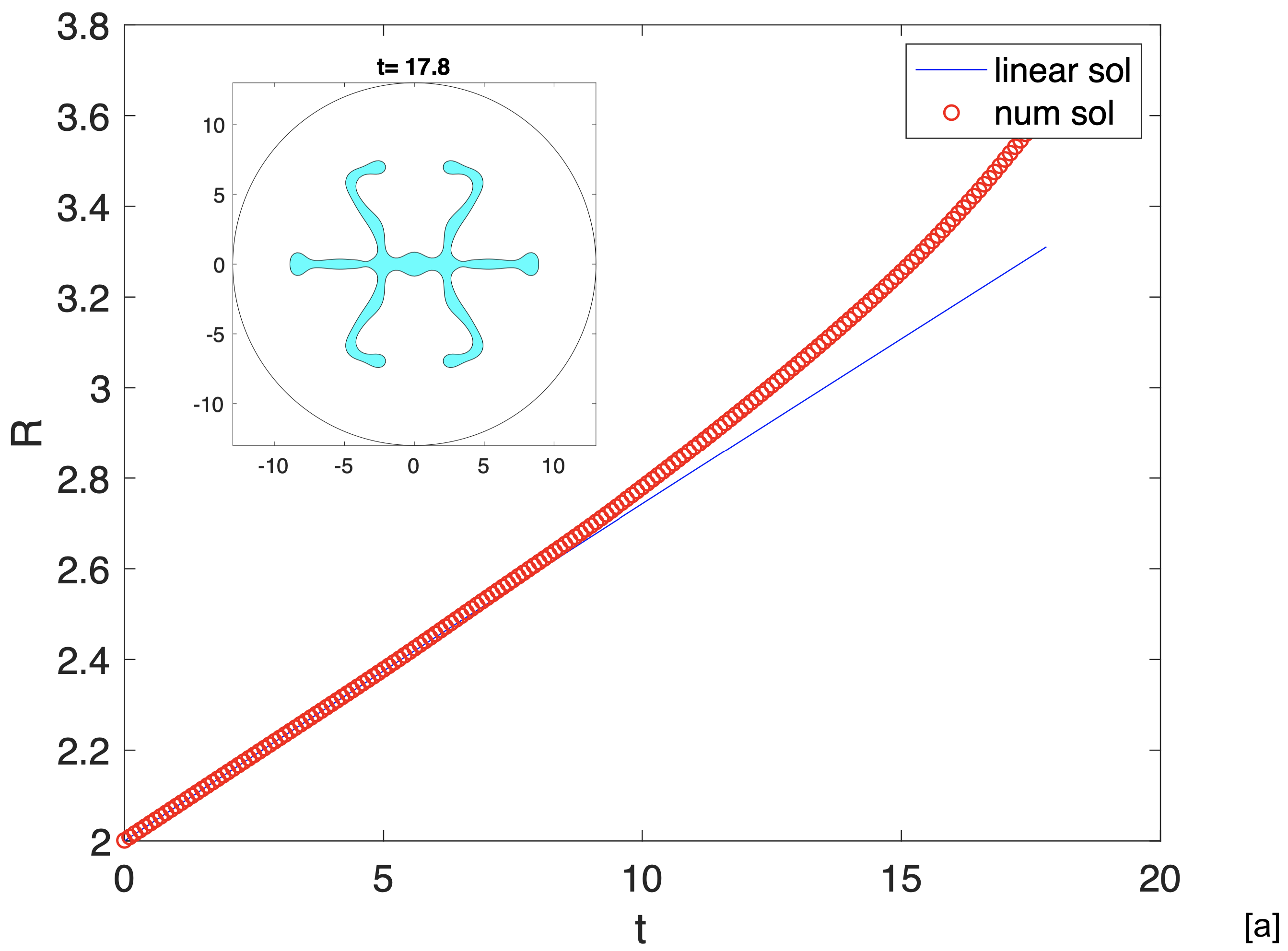}
\end{minipage}
\\
\begin{minipage}{1.0\linewidth}
\centering
\includegraphics[width=\textwidth]{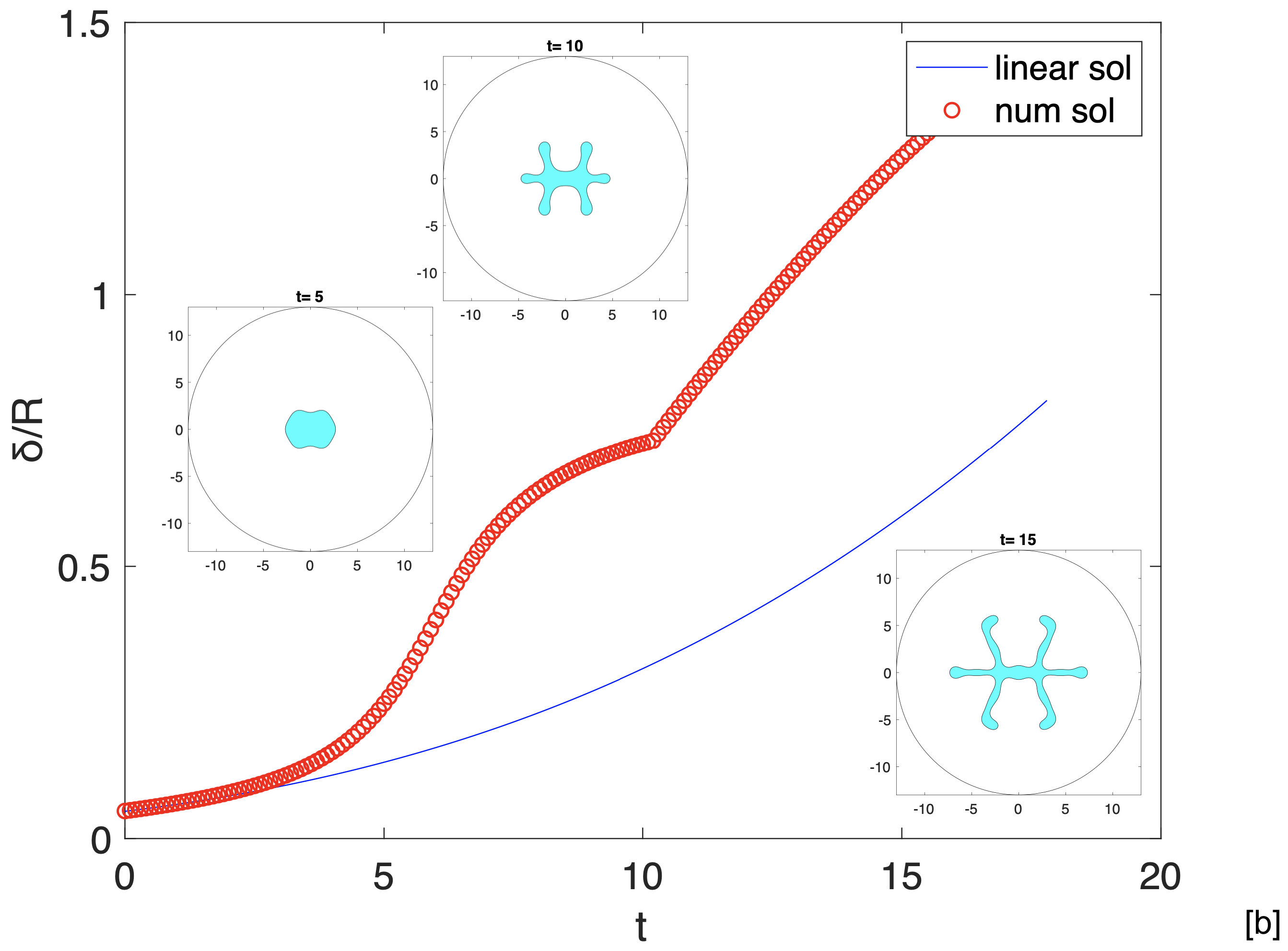}
\end{minipage}
\caption{A comparison between linear theory (blue curves) and the nonlinear simulations (red circles) for the effective radius $R$ in [a] and shape factor $\frac{\delta}{R}$ in [b]. Insets show the nonlinear tumor morphologies. The parameters are $D=1$, $\lambda=0.01$, $\chi_{\sigma}=5$, $\mathcal{P}=0.5$, $\mathcal{A}=0$,  $\mathcal{G}^{-1}=0.001$. The  far-field boundary is a circle $R_\infty=13$, and the initial tumor boundary is $r=2.0+ 0.1\cos (2 \theta)$. Here, $N=512$ and $\Delta t=0.005$.}
\label{fig:4}
\end{figure}

We next compare the nonlinear simulation with linear theory for longer times. The results are shown in Fig. \ref{fig:4} where we consider the case with the microenvironmental parameters $D=1$, $\lambda=0.01$, other parameters $\mathcal{P}=0.5$, $\mathcal{A}=0$, $\chi_{\sigma}=5$, $\mathcal{G}^{-1}=0.001$, simple (circular) far-field boundary $R_\infty=13$, and initial tumor surface $r=2.0+ 0.1\cos (2 \theta)$. While there is good agreement between the linear and nonlinear results at early times, both the effective radius and shape factors are under predicted by linear theory at later times. The nonlinear tumors, shown as insets, show the development of branched tubular structures similar to those observed in simulations \cite{cristini2009,garcke2016cahn} and experiments \cite{Pennacchietti2003}. In the next section, we will investigate how the far-field geometry influences tumor progression.

\subsection{Nonlinear Simulation with Non-Circular, Far-Field Geometries}
\label{sec:4.3}

In this section we demonstrate the dependence of tumor growth on the geometry of the far-field boundary under two prototypical microenvironmental conditions. In the first, which we call nutrient-poor, we take $D=1$ and $\lambda=0.01$, e.g.,  the ECM density in the host is similar to that in the tumor but the host cell density is smaller (since $\lambda$ is small). In the second, which we call nutrient-rich, we take $D=100$ and $\lambda=1$, e.g., the host ECM density is much smaller than that in the tumor but the host cell density is higher. In both cases, however, the diffusional penetration length is $\Lambda=10$. Unless otherwise specified, we take $\mathcal{P}=0.5$ and $\mathcal{G}^{-1}=0.001$.

\subsubsection{Nutrient-Poor Regime ($D=1$)}
\label{subsec:4.3.1}

\paragraph{Symmetric 5-fold Far-field Geometry.\\}

\begin{figure}
\includegraphics[width=\textwidth]{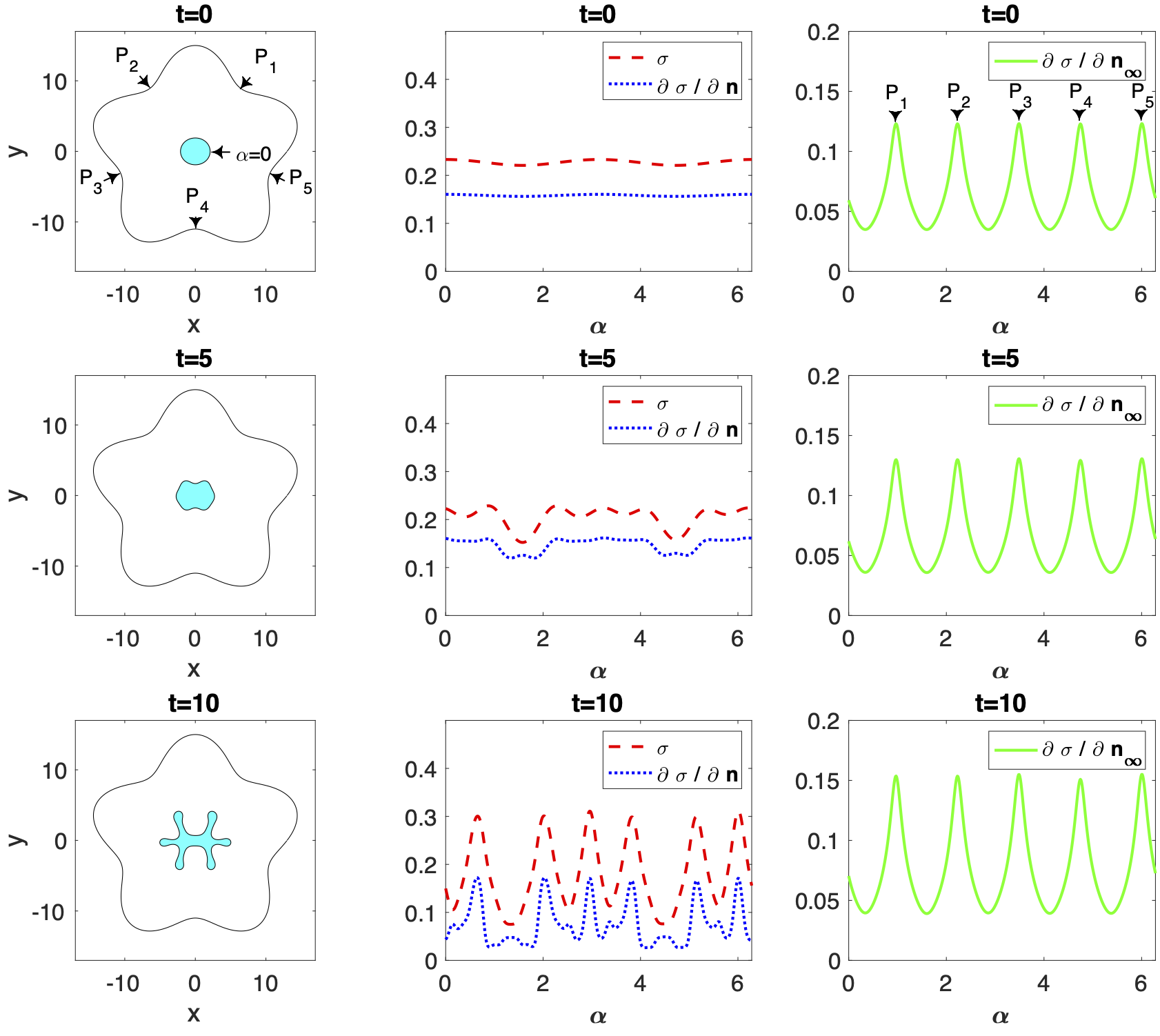}
\caption{Tumor morphologies (first column), the nutrient concentrations $\sigma$ and fluxes $\frac{\partial \sigma}{\partial \mathbf{n}}$ at $\Gamma$ (second column), and the far-field fluxes $\frac{\partial \sigma}{\partial \mathbf{n_\infty}}$ at $\Gamma_\infty$ (3rd column) in the nutrient-poor regime ($D=1$) and apoptosis $\mathcal{A}=0$.  The five-fold symmetric far-field boundary is $R_\infty=13+2\cos{(5\theta-\frac{\pi}{2})})$. The remaining parameters are $\lambda=0.01$, $\chi_{\sigma}=5$, $\mathcal{P}=0.5$,   and $\mathcal{G}^{-1}=0.001$. The initial tumor boundary is $r=2.0+ 0.1\cos (2 \theta)$. Here, $N=512$ and $\Delta t=0.005$.}
\label{fig:5}
\end{figure}

As a first test of the influence of a complex far-field geometry on tumor growth, we take a symmetric 5-fold perturbed circle instead of a circle as the far-field boundary with all other parameters as in Fig. \ref{fig:4}. We present the results in Fig. \ref{fig:5} where the first column shows tumor morphologies, the second column shows the nutrient concentrations $\sigma$ and fluxes $\frac{\partial\sigma}{\partial \mathbf{n}}$ at $\Gamma$, and the third column shows the fluxes $\frac{\partial\sigma}{\partial \mathbf{n}_\infty}$ at $\Gamma_\infty$. In the first column of Fig. \ref{fig:5} we see the tumor develops a similar morphology as in Fig. \ref{fig:4} with slight differences in the tubular structures. In the second column, we notice that both the nutrient concentration $\sigma$ and flux $\frac{\partial\sigma}{\partial\mathbf{n}}$ are high around the fingers, which indicates that in this case, the regions of highest proliferation and taxis are aligned. In the third column the nutrient fluxes $\frac{\partial\sigma}{\partial\mathbf{n}_\infty}$ at $\Gamma_\infty$ have similar profiles but the magnitudes of the peaks increase as time $t$ increases.
Next, we test different symmetric far-field geometries and different initial tumor shapes.

\paragraph{Symmetric 3,4,5,6-fold far-field geometries.\\}

\begin{figure}
\includegraphics[width=\textwidth]{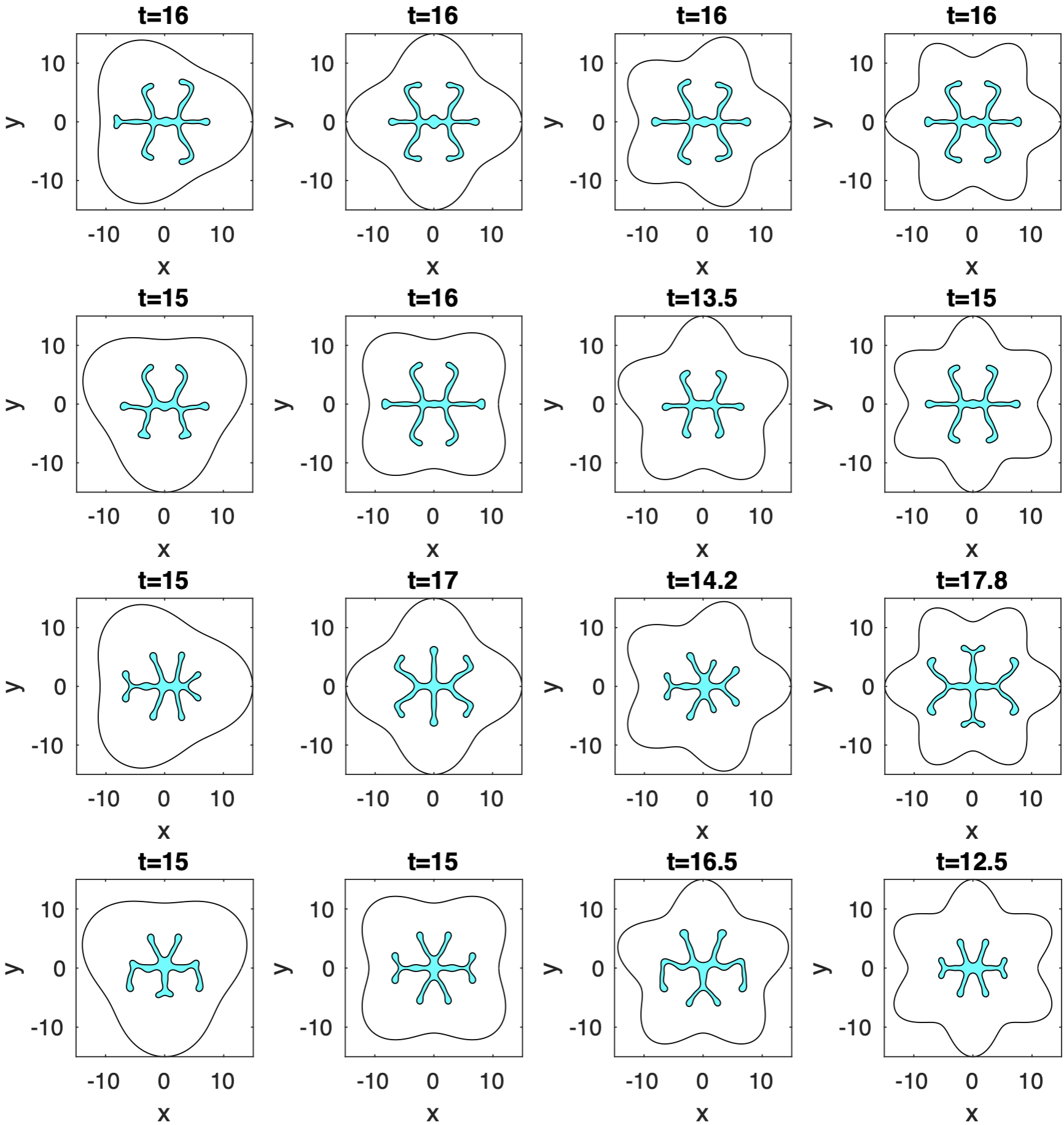}
\caption{Tumor morphologies with different symmetric far-field geometries in the nutrient-poor regime ($D=1$). The far-field boundaries are $R_\infty=13+2\cos(k\theta),k=3,4,5,6$ (Rows 1,3);
$R_\infty=13+2\cos(k\theta-\pi/2),k=3,5,~R_\infty=13+2\cos(k\theta-\pi),k=4,6$ (Rows 2,4). The initial tumor boundaries are $r=2.0+ 0.1\cos (2 \theta)$ (Rows 1,2) and $\frac{x^2}{2.1^2}+\frac{y^2}{1.9^2}=1$ (Rows 3,4). The remaining parameters are $\lambda=0.01$, $\chi_{\sigma}=5$, $\mathcal{P}=0.5$, $\mathcal{A}=0$ and $\mathcal{G}^{-1}=0.001$. Here, $N=512$, and $\Delta t=0.005$.}
\label{fig:6}
\end{figure}

In Fig. \ref{fig:6}, we present tumor morphologies at similar sizes under the same growth conditions but using  different symmetric far-field boundaries: $R_\infty=13+2\cos(k\theta),k=3,4,5,6$ (Row 1,3) and $R_\infty=13+2\cos(k\theta-\pi/2),k=3,5,R_\infty=13+2\cos(k\theta-\pi),k=4,6$ (Row 2,4). The parameters are the same as in Fig. \ref{fig:5}, where in rows 1 and 2 the initial tumor boundary is the perturbed circle $r=2.0+ 0.1\cos (2 \theta)$ while in rows 3 and 4 the initial tumor boundary is the ellipse $\frac{x^2}{2.1^2}+\frac{y^2}{1.9^2}=1$. Hence, by comparing the evolution within each row, and between rows 1 and 2 and rows 3 and 4, we can see the effect of the far-field boundary shapes. By comparing rows 1 and  3 and rows 2 and 4, we can see the effect of the different initial shapes.

All the tumors develop tubular structures. When the initial condition is the perturbed circle (rows 1 and 2), the far-field geometry has limited influence on tumor morphologies consistent with that observed in Figs. \ref{fig:4} and \ref{fig:5}. However, when the initial condition is an ellipse, which contains many modes, the morphologies are much more sensitive to the far-field geometry.  As a consequence, the morphological instability drives growth of more modes than in the perturbed circular case and this growth is more strongly coupled with the far-field geometry. We note that the tumor morphology for the 4-fold symmetric far-field geometry in row 4 (elliptic initial condition) is similar to that observed in Fig. 11(b) in \cite{cristini2009} where the far-field geometry is a square. When the chemotaxis coefficient $\chi_\sigma$ is increased, the tubular structures develop faster and are thinner but the overall morphologies are similar (see "Appendix E").

\paragraph{Non-symmetric complex far-field geometry.\\}

\begin{figure}
\includegraphics[width=\textwidth]{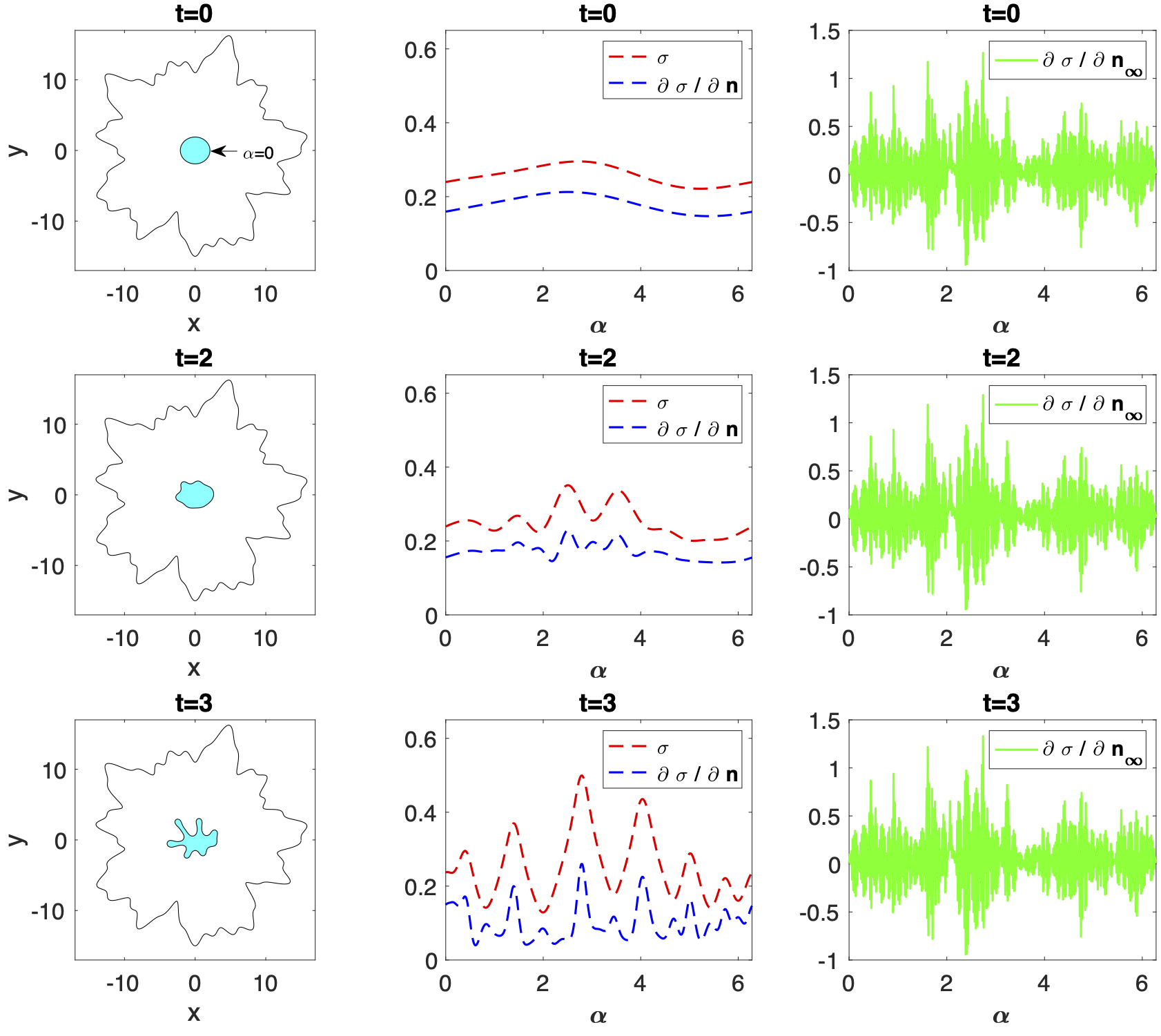}
\caption{Tumor morphologies (first column), the nutrient concentrations $\sigma$ and fluxes $\frac{\partial \sigma}{\partial \mathbf{n}}$ at $\Gamma$ (second column), and the far-field fluxes $\frac{\partial \sigma}{\partial \mathbf{n_\infty}}$ at $\Gamma_\infty$ (third column) in the nutrient-poor regime ($D=1$) and apoptosis $\mathcal{A}=0$.  The non-symmetric far-field boundary is $R_\infty=13+2(\cos^3{(5\theta)}+\sin^3{(11\theta)})$. The remaining parameters are $\lambda=0.01$, $\chi_{\sigma}=10$, $\mathcal{P}=0.5$,   and $\mathcal{G}^{-1}=0.001$. The initial tumor boundary is $\frac{x^2}{2.1^2}+\frac{y^2}{1.9^2}=1$. Here, $N=512$, and $\Delta t=0.005$.}
\label{fig:7}
\end{figure}

In Fig. \ref{fig:7} we show the tumor evolution in the non-symmetric complex far-field boundary $R_\infty=13+2(\cos^3{(5\theta)}+\sin^3{(11\theta)})$. The initial tumor shape is an ellipse and the parameters are the same as in Fig. \ref{fig:6} except that $\chi_\sigma=10$. The first column shows the tumor morphological development, the second shows the detailed nutrient and flux distributions on $\Gamma$ and the third shows the nutrient flux $\partial\sigma/\partial \mathbf{n}_\infty$ at $\Gamma_\infty$. The spatial distribution of $\partial\sigma/\partial \mathbf{n}_\infty$ is spiky. It appears to be very noisy, with small scale variations driven by the complexity of the far-field boundary, but the distribution of $\partial\sigma/\partial \mathbf{n}_\infty$ is in fact deterministic, resolved by the computational mesh, and changes only slowly in time. Compared to Fig. \ref{fig:5}, where a circular far-field boundary was used, the magnitude of the flux $\partial\sigma/\partial \mathbf{n}_\infty$ is about 10 times larger here.

As in the previous cases, the growing tumor develops tubular structures although the detailed tumor morphology is asymmetric and influenced by the far-field geometry.  Compared to the cases shown in Fig. \ref{fig:6} and in Appendix E, the curvature of the tumor/host interface is increased due to long-range interactions with the far-field boundary through the nutrient gradients, even though the tumor is still far from the outer boundary. This requires the use of smaller time and space steps to continue solving the system beyond the time shown. As seen in the second column of  Fig. \ref{fig:7}, both the nutrient concentration $\sigma$ and flux $\frac{\partial\sigma}{\partial\mathbf{n}}$ are high around these tips of the tubular structures, as seen earlier in Fig. \ref{fig:5}, although the magnitudes here are larger.

\subsubsection{Nutrient-Rich Regime ($D=100$)\\}
\label{subsec:4.3.2}

Next, we consider tumor evolution in the nutrient-rich regime with $D=100$. We take $\lambda=1$ to maintain the same diffusional penetration length $\Lambda$. We also take $\mathscr{P}=0.5$ and $\mathscr{G}^{-1}=0.001$ as before.

\paragraph{Compact growth with $\mathcal{A}=0$\\}

\begin{figure}
\includegraphics[width=\textwidth]{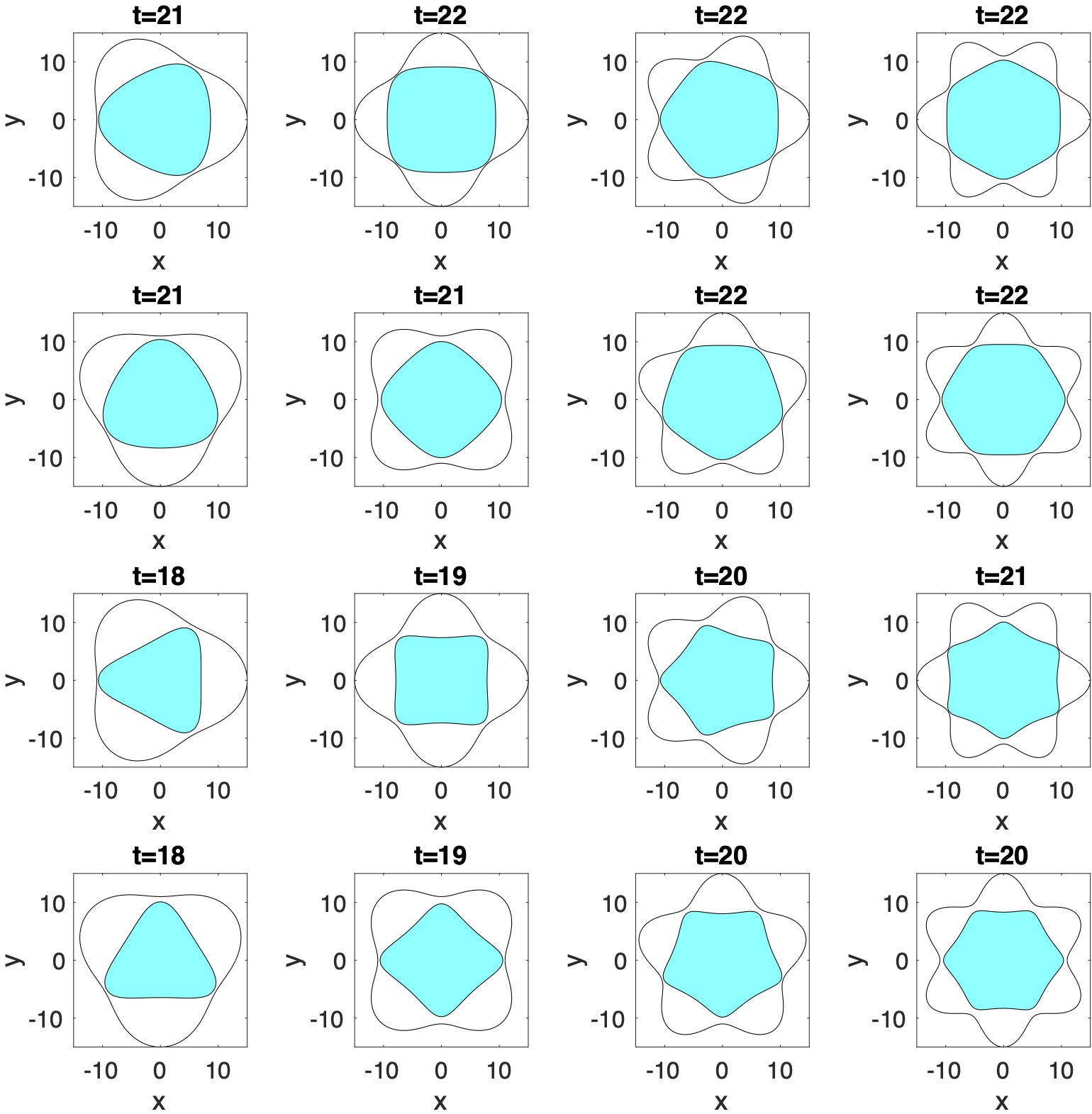}
\caption{Tumor morphologies with different symmetric far-field geometries in the nutrient-rich regime ($D=100$). The far-field boundaries are $R_\infty=13+2\cos(k\theta),k=3,4,5,6$ (Rows 1,3);
$R_\infty=13+2\cos(k\theta-\pi/2),k=3,5,~R_\infty=13+2\cos(k\theta-\pi),k=4,6$ (Rows 2,4). In rows 1 and 2, $\chi_{\sigma}=5$ and in rows 3 and 4, $\chi_{\sigma}=10$. The initial tumor shape is the perturbed circle $r=2.0+0.1 \cos (2 \theta)$. The remaining parameters are $\lambda=1$, $\mathcal{P}=0.5$, $\mathcal{A}=0$,  $\mathcal{G}^{-1}=0.001$. Here, $N=512$, and $\Delta t=0.005$.
}
\label{A zero}
\end{figure}

In Fig. \ref{A zero}, we present tumor morphologies for different symmetric, far-field boundaries similar to Fig. \ref{fig:6}, where $D=1$. Here, the initial tumor shape is the perturbed circle $r=2.0+0.1 \cos (2 \theta)$. Rows 1 and 2 have $\chi_\sigma=5$ while rows 3 and 4 use  $\chi_\sigma=10$. In contrast to the $D=1$ case shown in Fig. \ref{fig:6} and in "Appendix E," when $D=100$ the tumors grow faster, develop compact morphologies, and the shapes are strongly dependent on the far-field geometry. Figure  \ref{A zero} also shows that $\chi_\sigma$ has a limited influence on the tumor morphologies because the larger $D$ reduces nutrient gradients making more nutrient available to the tumor. The tumors with larger $\chi_\sigma$ (rows 3 and 4) are slightly more unstable and develop regions of negative curvature.

\begin{figure}
\begin{minipage}{1.0\linewidth}
\centering
\includegraphics[width=\textwidth]{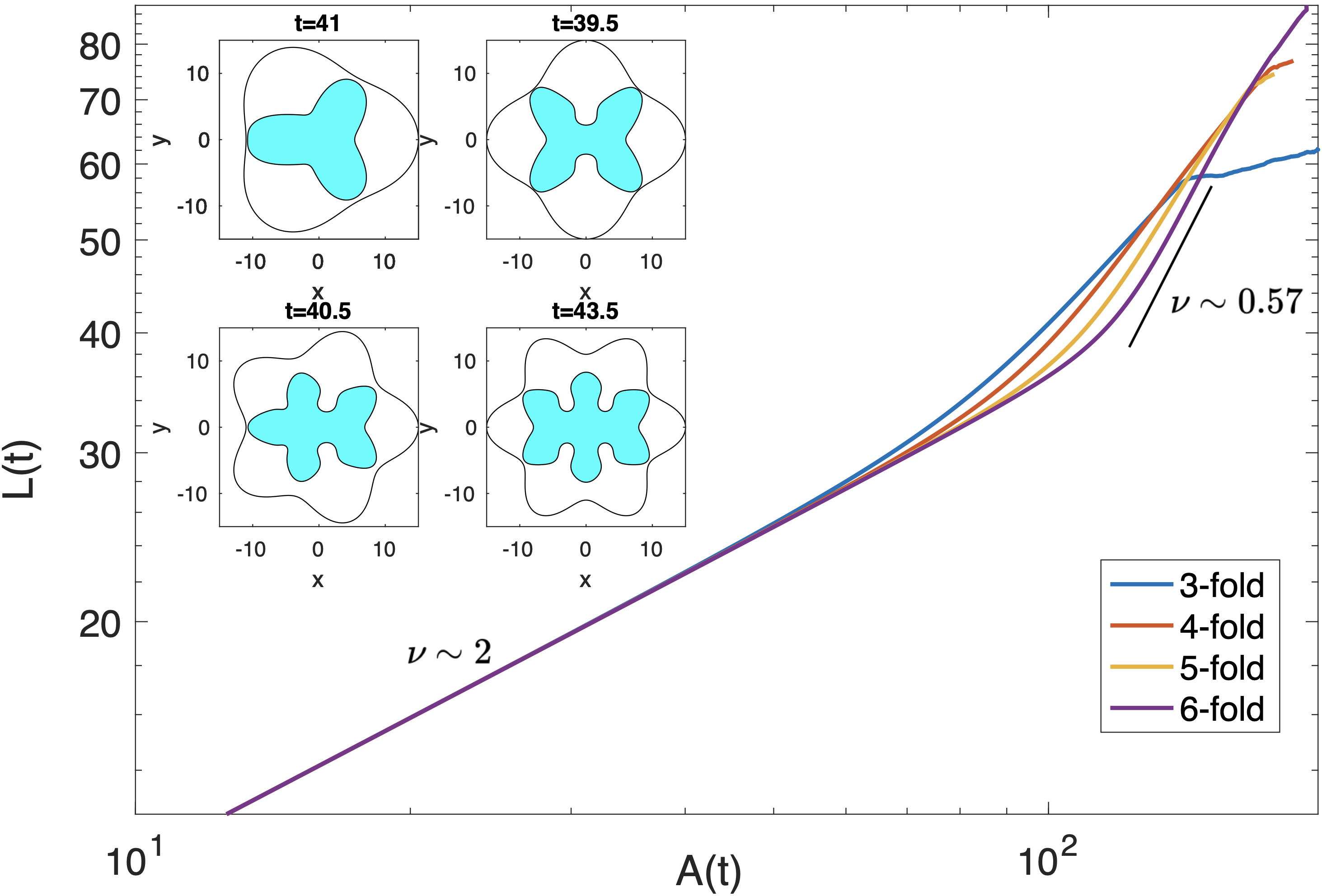}[a]
\end{minipage}

\begin{minipage}{1.0\linewidth}
\centering
\includegraphics[width=\textwidth]{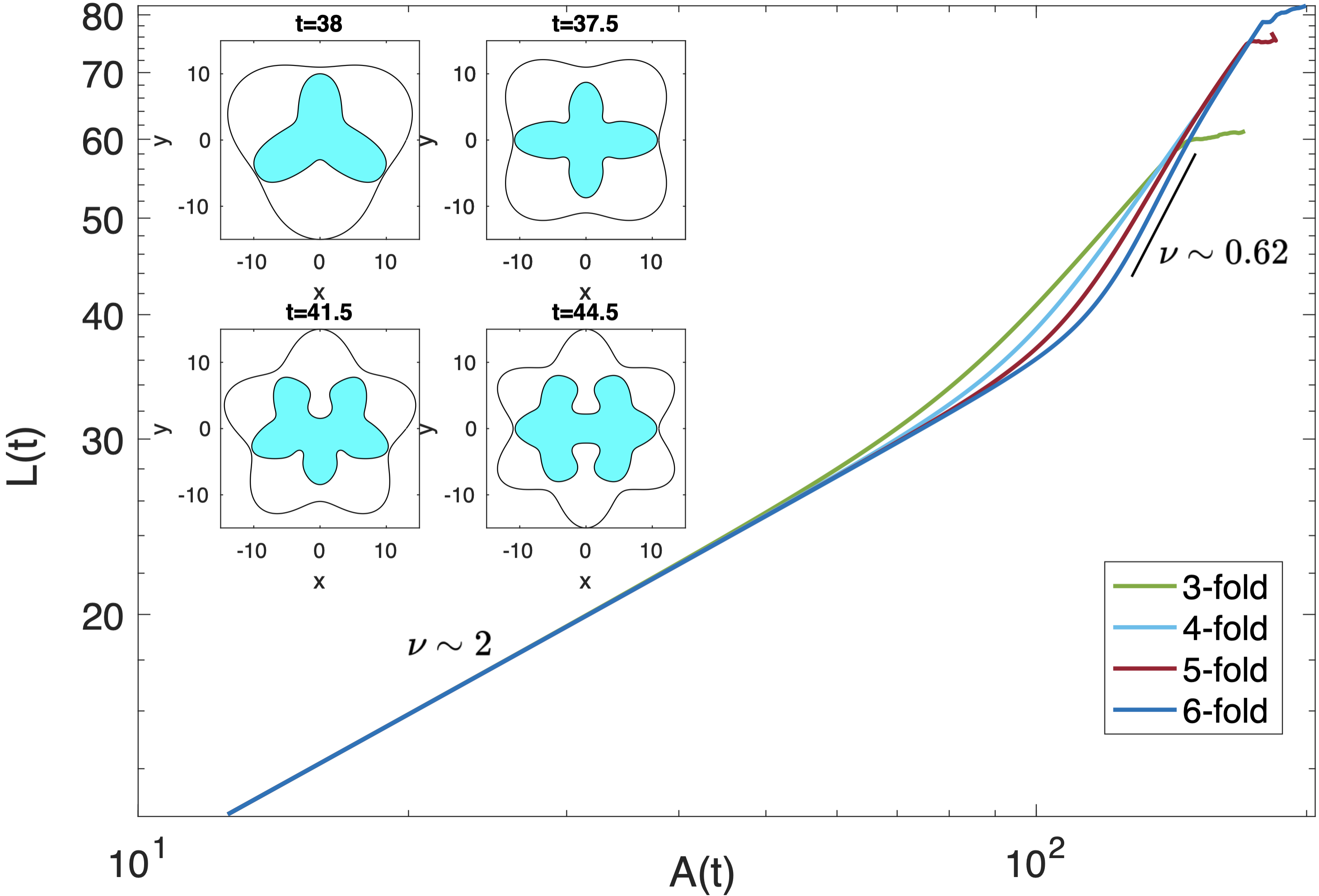}[b]
\end{minipage}

\caption{Tumor morphologies (insets) and the relation between the tumor arclengths and areas during the evolution under different symmetric far-field boundaries in the nutrient-rich regime ($D=100$) and apoptosis $\mathcal{A}=0.1$. The far-field boundaries are  [a] $R_\infty=13+2\cos(k\theta),k=3,4,5,6;$ [b] $R_\infty=13+2\cos(k\theta-\pi/2),k=3,5,R_\infty=13+2\cos(k\theta-\pi),k=4,6$. The initial tumor surface is the perturbed circle $r=2.0+ 0.1\cos (2 \theta)$. The remaining parameters are $\lambda=1$, $\chi_{\sigma}=5$, $\mathcal{P}=0.5$,  $\mathcal{G}^{-1}=0.001$. Here, $N=512$, and $\Delta t=0.005$.}
\label{fig:Area-arclength}
\end{figure}

\paragraph{Unstable growth with $\mathcal{A}=0.1$.\\}
When the apoptosis rate $\mathcal{A}$ is increased, the instability is enhanced as predicted by linear theory. This leads to slower growth and the development of tumor boundaries with increased negative curvature, as seen in the insets of Fig. \ref{fig:Area-arclength}. In Fig. \ref{fig:Area-arclength}, $\mathcal{A}=0.1$ and $\chi_\sigma=5$. The initial condition is the perturbed circle $r=2.0+ 0.1\cos (2 \theta)$ and all the other parameters are the same as in Fig. \ref{A zero}. The insets in Fig. \ref{fig:Area-arclength} show that tumors are strongly influenced by the far-field boundary geometry and develop buds that grow into thick, invasive fingers that protrude into the surrounding host tissue. The larger the tumor grows, the more the arclength $L$ to area enclosed $A$ ratio of the tumor boundary increases. At early times, the growth is compact with $A\approx L^2$. but as the thick fingers develop, the length increases more rapidly with a relation $A\approx L^\nu$ with $\nu\approx 0.6$, as labeled. In all four cases, the growth dynamics is similar, which suggests there may be a universal scaling between the area and the arclength.

\paragraph{Unstable growth with $\mathcal{A}=0.2$.\\}

\begin{figure}
\includegraphics[width=\textwidth]{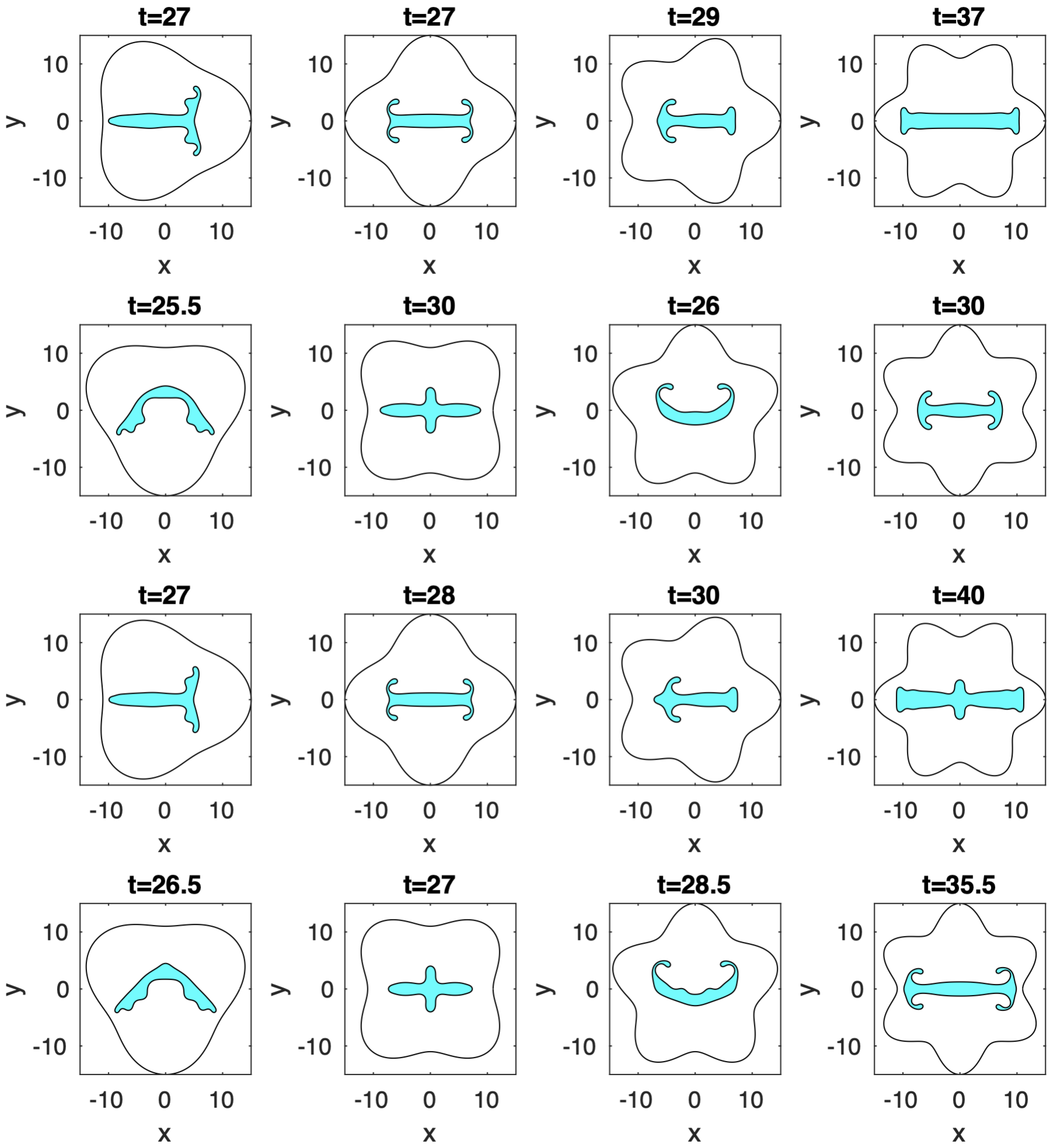}
\caption{Tumor morphologies with different symmetric far-field geometries in the nutrient-rich regime ($D=100$) and apoptosis $\mathcal{A}=0.2$. The far-field boundaries are $R_\infty=13+2\cos(k\theta),k=3,4,5,6$ (Rows 1,3);
$R_\infty=13+2\cos(k\theta-\pi/2),k=3,5,~R_\infty=13+2\cos(k\theta-\pi),k=4,6$ (Rows 2,4). The initial tumor boundaries are $r=2.0+ 0.1\cos (2 \theta)$ (Rows 1,2), and $\frac{x^2}{2.1^2}+\frac{y^2}{1.9^2}=1$ (Rows 3,4). The remaining parameters are $\lambda=1$, $\chi_{\sigma}=10$, $\mathcal{P}=0.5$ and $\mathcal{G}^{-1}=0.001$. Here, $N=512$, and $\Delta t=0.005$.}
\label{fig:D100A02-3456farfield-c10}
\end{figure}

Next we increase apoptosis to $\mathcal{A}=0.2$. In Fig. \ref{fig:D100A02-3456farfield-c10}, we present tumor morphologies for different symmetric, far-field geometries and different initial conditions using the same arrangement in Fig. \ref{fig:6} where $D=1$. That is, the initial tumor shapes are perturbed circles in rows 1 and 2 and ellipses in rows 3 and 4. In all cases, $\chi_\sigma=10$. As seen in Fig. \ref{fig:D100A02-3456farfield-c10}, apoptosis generates strong instabilities but the shapes produce less branched compared to the nutrient-poor case ($D=1$). The tumor morphologies are qualitatively similar to those obtained in the nutrient-rich regime considered in \cite{macklin2007nonlinear} where chemotaxis was absent and the far-field boundary was a circle. In contrast here, the tumor shapes are seen to be sensitively dependent on the far-field geometry. Due to the large $D$, the tumor morphologies are much less sensitive to the initial shape, however.

\paragraph{Complex, non-symmetric far-field geometry with $\mathcal{A}=0.2$.\\}
\begin{figure}
\includegraphics[width=\textwidth]{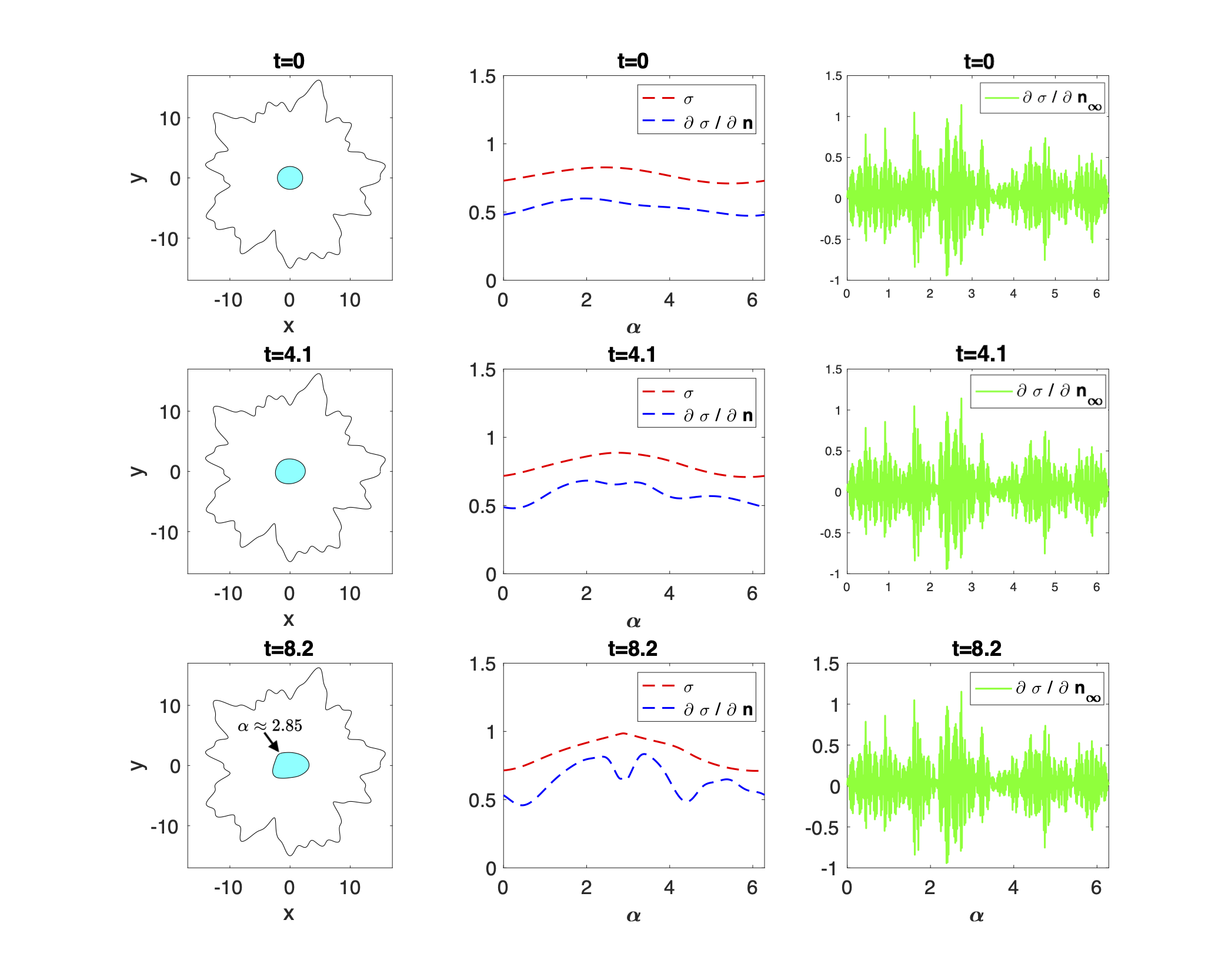}
\caption{Tumor morphologies (first column), the nutrient concentrations $\sigma$ and fluxes $\frac{\partial \sigma}{\partial \mathbf{n}}$ at $\Gamma$ (second column), and the far-field fluxes $\frac{\partial \sigma}{\partial \mathbf{n_\infty}}$ at $\Gamma_\infty$ (third column) in the nutrient-rich regime ($D=100$) with apoptosis $\mathcal{A}=0.2$. The non-symmetric far-field boundary is $R_\infty=13+2(\cos^3{(5\theta)}+\sin^3{(11\theta)})$. The remaining parameters are $\lambda=1$, $\chi_{\sigma}=10$, $\mathcal{P}=0.5$ and $\mathcal{G}^{-1}=0.001$. The initial tumor boundary is $\frac{x^2}{2.1^2}+\frac{y^2}{1.9^2}=1$. Here, $N=512$, and $\Delta t=0.005$. The arrow in row 3, column 1 labels the position on the tumor boundary where a corner-like shape with high curvature develops.}
\label{fig:D100A02unsym}
\end{figure}

In Fig. \ref{fig:D100A02unsym}, we present a tumor evolving in the complex, non-symmetric far-field geometry considered in Fig. \ref{fig:7} in the nutrient-poor regime with $D=1$. Here, using $D=100$, we find the tumor grows compactly but develops a corner-like shape with high curvature at $t\approx 8.2$ near $\alpha\approx 2.85$, as labeled.  As seen in the second column, the nutrient $\sigma$ and the flux $\partial\sigma/\partial \mathbf{n}_\infty$ are out of phase with the peak of $\sigma\approx 0.98$ occurring at a local minima of $\partial\sigma/\partial \mathbf{n}_\infty$. This gives rise to a competition between proliferation and chemotaxis, which drives the formation of the near corner. Smaller spatial and temporal grids are needed to determine whether a singularity develops at a finite time. This will be considered in future work. The behavior of the flux at the far-field boundary (3rd column) is similar to that shown in Fig. \ref{fig:7}.

\subsubsection{Comparison Between Nutrient-Poor $(D=1)$ and Nutrient-Rich $(D=100)$ Regime}
In general, tumors in the nutrient-poor regime \((D=1)\)  are much less-dependent
 on the boundary geometries than tumors in the nutrient-rich regime \((D=100)\). This is because the nutrient decays more rapidly away from the far-field domain boundary in the nutrient-poor regime (e.g., Fig. \ref{fig:7} where $D=1$) compared to the nutrient-rich regime (e.g., Fig. \ref{fig:D100A02unsym} where $D=100$) and hence the nutrient distribution on the tumor boundary is less dependent on the far-field geometry than in the nutrient-rich regime. Therefore, in the nutrient-poor regime the tumor morphology is
 more controlled by the internal parameters--chemotaxis $\chi_\sigma$, apoptosis $\mathcal{A}$-- than the far-field geometry. Both $\chi_\sigma$ and $\mathcal{A}$ promote morphological instabilities. For example, as  $\mathcal{A}$ is increased, the morphologies in the nutrient-rich regime range from being entirely determined by the boundary geometry (e.g., Fig. \ref{A zero} where \(D=100,\mathcal{A}=0\)) to rod-like shapes that are much less sensitive to the boundary (e.g., Fig. \ref{fig:D100A02-3456farfield-c10} where \(D=100,\mathcal{A}=0.2\)) because shape changes due to cell death $\mathcal{A}>0$ compete with those induced by the boundary geometries. In addition, the effect of chemotaxis in the nutrient-poor regime generates larger instability for the same value of $\chi_\sigma$ than in the nutrient-rich regime because the nutrient gradients are larger in the nutrient-poor regime.

\section{Conclusion}
\label{sec:6}
In this paper, we have developed, analyzed and solved numerically a tumor growth model that accounts for complex far-field geometries, which mimic heterogeneous vascular distributions. The model incorporates cell proliferation, death and chemotaxis up gradients of nutrients that are transported diffusionally from the far-field boundary and uptaken by both tumor and host cells. Linear analysis, though limited to a simple far-field geometry (performed here for a circular far-field boundary), reveals the presence of rich pattern formation mechanisms via unstable tumor growth.

To gain insight into the nonlinear solutions, we developed a novel boundary integral method to accurately and efficiently simulate the system. A direct layer potential representation was used for two-phase nutrient equation, which enables us to obtain the value of the nutrient concentration and the nutrient flux on the interface and on the far-field boundary accurately by solving three coupled integral equations. The tumor interface was evolved using a semi-implicit time-stepping method developed previously (e.g., \cite{hou1994removing,cristini2003nonlinear}). The method is spectrally accurate in space and second order accurate in time.

With the advantage of boundary integral methods in addressing the complex tumor and far-field  geometries, our nonlinear simulations explored various unstable morphologies caused by vascular heterogeneities, nutrient diffusion constants, and cancer cell apoptosis. When the nutrient diffusion constants are comparably low both inside and outside the tumor, the tumor is less sensitive to the far-field geometry, and the tumor may grow to a thin tubular star shape, similar to those found in previous theoretical \cite{cristini2009,garcke2016cahn} and experimental \cite{Pennacchietti2003} studies. When the nutrient diffusion constant in the host tissue is far larger than that in the tumor, the tumor morphology is more sensitive to the far-field geometry, and the tumor would grow to a bulky compact shape if the cancer cell apoptosis rate is low, but to a rod-like shape if the apoptosis rate is high.

In future work, we can investigate the dynamics using other penetration lengths $\Lambda$ together with other components of the microenvironment such as stromal and immune cells. In addition, we can consider passive and active cell movement in the host tissue, including the motion of the far-field boundary (e.g., two moving boundaries in the system). Further, the nutrient concentration on the far-field boundary need not be constant and the cell division and uptake rates need not be uniform, as assumed here. Another immediate extension is to use the Stokes equations to model the dynamics of the tumor and host tissue. This will highlight the effects of viscosity contrast (e.g., tissue stiffness) between the tumor and the host tissues. The regulation of cell fates and motility, proliferation and apoptosis rates by mechanical and thermal stresses can also be incorporated. Finally, while we presented the results in two dimensions, similar behaviors are expected to hold qualitatively in three dimensions, and we plan to perform full 3D simulations to confirm this.

\begin{acknowledgements}
We would like to acknowledge the referee's insightful suggestions and contributions.
S. L. acknowledges the support from the National Science Foundation, Division of Mathematical Sciences Grant DMS-1720420. S. L. was also partially supported by Grant ECCS-1307625. M. L. acknowledges the F. R. Buck McMorris Summer Research support from the College of Science, IIT. C. L. is partially supported by the National Science Foundation, Division of Mathematical Sciences Grant DMS-1759536. J.L. acknowledges partial support from the NSF through Grants DMS-1714973, DMS-1719960, and DMS-1763272 and the Simons Foundation (594598QN) for a NSF-Simons Center for Multiscale Cell Fate Research. J.L. also thanks the National Institutes of Health for partial support through Grants 1U54CA217378-01A1 for a National Center in Cancer Systems Biology at UC Irvine and P30CA062203 for the Chao Family Comprehensive Cancer Center at UC Irvine.

\end{acknowledgements}
\appendix
\renewcommand{\thesection}{\Alph{section}}
\renewcommand{\thesubsection}{A.\arabic{subsection}}
\section*{Appendix A: Linear Stability Analysis}
\label{appendix}
The governing equations are

\begin{align}
\Delta \sigma-\mu_{\mathrm{i}}^{2} \sigma
&=0 \text { in } \Omega_{\mathrm{i}},\\
\ [\sigma]_{\Gamma}
&=0\label{noslipperybc},\\
\left[D \frac{\partial \sigma}{\partial n}\right]_{\Gamma}
&=0\label{fluxbc},\\
\sigma|_{\Gamma_{\infty}}
&=1.\label{farfieldbc}
\end{align}
and
\begin{align}{}
\Delta p
&=0 \text { in } \Omega_{1},\\
{p}
&=\mathcal{G}^{-1} \kappa+\left(\mathcal{P}-\chi_{\sigma}\right) \sigma-\mathcal{A} \frac{\mathbf{x} \cdot \mathbf{x}}{2 {d}}\text{ on } \Gamma.\label{pressurebc}
\end{align}
\begin{equation}
V=-\frac{\partial p}{\partial n}+\mathcal{P}\frac{\partial \sigma}{\partial n}-\mathcal{A} \frac{\mathbf{n} \cdot\mathbf{x}}{d}\text{ on } \Gamma.\label{normalvel}
\end{equation}Consider a perturbed tumor interface $\Gamma$:
\begin{equation}
r(t)=R(t)+\delta(t) e^{i l \theta}\label{perturbedcircle}.
\end{equation}
In cylindrical coordinates, the modified Helmholtz equation is
\begin{equation}
r^{-1}\left(r \sigma_{r}\right)_{r}+r^{-2} \sigma_{\theta \theta}+r^{-2} \sigma_{z z}-\mu_{i}^{2} \sigma=0 \text { in } \Omega_{i}.
\end{equation}Assuming axial symmetry, $\textit{e.g.}$, $\sigma=\sigma(r, \theta)$ is independent of $z$, then

\begin{equation}
r^{-1}\left(r \sigma_{r}\right)_{r}+r^{-2} \sigma_{\theta \theta}-\mu_{i}^{2} \sigma=0\label{MHHeq}.
\end{equation}
\subsection{Radial Solutions}
We first consider the radial solution, \textit{i.e.}, $\sigma=\sigma(r)$, then $\eqref{MHHeq}$ reduces to the modified Bessel differential equation
\begin{equation}
\left(r^{2} \frac{d^{2}}{d r^{2}}+r \frac{d}{d r}-\mu_{i}^{2} r^{2}\right) \sigma(r)=0 \text { in } \Omega_{i}.
\end{equation}
Recall the general form of modified Bessel differential equation is
\begin{equation}
\left(x^{2} \frac{d^{2}}{d x^{2}}+x \frac{d}{d x}-\left(\beta^{2} x^{2}+n^{2}\right)\right) y(x)=0.
\end{equation}
The general solutions are
\begin{eqnarray}{}
y&=& a_{1} J_{n}(-i \beta x)+a_{2} Y_{n}(-i \beta x) \nonumber\\
&=&c_{1} I_{n}(\beta x)+c_{2} K_{n}(\beta x),
\label{bessels}
\end{eqnarray}
where $J_{n}(x)$ is the Bessel function of the first kind, $Y_{n}(x)$ is the Bessel function of the second kind, $I_{n}(x)$ is a modified Bessel function of the first kind and $K_{n}(x)$ is a modified Bessel function of the second kind.

The following recurrence relations are useful in the linear analysis
\begin{eqnarray}{}
I_{n}^{\prime}(x)&=&\frac{1}{2}\left(I_{n-1}(x)+I_{n+1}(x)\right),\nonumber\\
I_{n}^{\prime}(x)&=&I_{n-1}(x)-\frac{n}{x} I_{n}(x)=\frac{n}{x} I_{n}(x)+I_{n+1}(x),\nonumber\\
I_{0}^{\prime}(x)&=&I_{1}(x).\label{besseliid}\\
K_{n}^{\prime}(x)&=&-\frac{1}{2}\left(K_{n-1}(x)+K_{n+1}(x)\right),\nonumber\\
K_{n}^{\prime}(x)&=&-K_{n-1}(x)-\frac{n}{x} K_{n}(x)=\frac{n}{x} K_{n}(x)-K_{n+1}(x),\nonumber\\
K_{0}^{\prime}(x)&=&-K_{1}(x).
\end{eqnarray}
The modified Bessel functions of the second kind all have the property that
\[K_{n}(x) \rightarrow \infty \text { as } x \rightarrow 0\nonumber.\]Hence, the constant $c_2$ in Eq. (\ref{bessels}) must be zero. We therefore obtain
\begin{align}
\sigma&=A_{1} I_{0}\left(\mu_{1} r\right) &\operatorname{in}\  \Omega_{1},\\
\sigma&=A_{2} I_{0}\left(\mu_{2} r\right)+A_{3} K_{0}\left(\mu_{2} r\right) &\operatorname{in}\  \Omega_{2}.
\end{align}Applying the boundary conditions $\eqref{noslipperybc},\eqref{fluxbc},\eqref{farfieldbc}$ on the circles $r=R,R_\infty$, we have
\begin{align}{}
A_{1} I_{0}\left(\mu_{1} R\right)&=A_{2} I_{0}\left(\mu_{2} R\right)+A_{3} K_{0}\left(\mu_{2} R\right),\\
1&=A_{2} I_{0}\left(\mu_{2} R_{\infty}\right)+A_{3} K_{0}\left(\mu_{2} R_{\infty}\right),\\
A_{1} \mu_{1} I_{1}\left(\mu_{1} R\right)&={D}\left(A_{2} \mu_{2} I_{1}\left(\mu_{2} R\right)-A_{3} \mu_{2} K_{1}\left(\mu_{2} R\right)\right).
\end{align}Solving for $A_{1}, A_{2}, A_{3}$, we obtain
\small
\begin{eqnarray}{}
A_1&=&
\frac{D}{(D \mu_{2} R I_{0}(R \mu_{1})
\left(
I_{1}(R \mu_{2})K_{0}(R_{\infty} \mu_{2})+K_{1}(R \mu_{2}) I_{0}(R_{\infty}\mu_{2})
\right)}\nonumber\\
&&+\mu_{1} R I_{1}
(R \mu_{1})
(K_{0}(R \mu_{2}) I_{0}
(R_{\infty} \mu_{2})
-I_{0}(R \mu_{2}) K_{0}
(R_{\infty} \mu_{2}))),\nonumber\\
A_2&=&
\frac{D \mu_{2} I_{0}(R \mu_{1})
K_{1}(R \mu_{2})
+\mu_{1} I_{1}(R \mu_{1}) K_{0}(R \mu_{2})}
{(D \mu_{2} I_{0} (R \mu_{1})
(I_{1} (R \mu_{2}) K_{0} (R_{\infty} \mu_{2})
+K_{1}(R \mu_{2}) I_{0} (R_{\infty} \mu_{2}) )
+\mu_{1} I_{1}(R \mu_{1})
(K_{0} (R \mu_{2}) I_{0} (R_{\infty} \mu_{2})
-I_{0}(R \mu_{2}) K_{0}(R_{\infty} \mu_{2}) ) )},\nonumber\\
A_3&=&
\frac{\mu_{1} I_{0} (R \mu_{2}) I_{1} (R \mu_{1})
-D \mu_{2} I_{0}(R \mu_{1}) I_{1}(R \mu_{2})}
{\mu_{1} I_{1}(R \mu_{1})((I_{0}(R \mu_{2}) K_{0}(R_{\infty} \mu_{2})-K_{0}(R \mu_{2}) I_{0}(R_{\infty} \mu_{2}))-D \mu_{2} I_{0}(R \mu_{1})(I_{1}(R \mu_{2}) K_{0}(R_{\infty} \mu_{2})+K_{1}(R \mu_{2}) I_{0}(R_{\infty}, \mu_{2}))}.\nonumber
\end{eqnarray}
\subsection{Perturbation of Radial Solutions}
Now we seek a solution of the modified Helmholtz equation on the perturbed circle given by $\eqref{perturbedcircle}$. Since $\delta$ is the perturbation size, following \cite{mullins1963morphological} we consider the Fourier expansion of the solution up to the first order in $\delta$:
\begin{equation}
\sigma(r, \theta)=\sigma_{i, 0}(r)+\delta e^{i l \theta} \sigma_{i, 1}(r) \text { in } \Omega_{i}.
\end{equation}Note here that $r, \theta$ and $\delta$ are all functions of time $t,$ \textit{i.e.,} $r=r(t), \theta=\theta(t), \delta=\delta(t)$.
Multiplying Eq. $\eqref{MHHeq}$ by $r^2$, we obtain
\begin{align}
\left(r^{2} \partial_{r}^{2}+r \partial_{r}+\partial_{\theta}^{2}-\mu_{i}^{2} r^{2}\right)\left(\sigma_{i, 0}(r)+\delta e^{i l \theta} \sigma_{i, 1}(r)\right)&=0 &\text { in } \Omega_{i},\\
\left(r^{2} \frac{d^{2}}{d r^{2}}+r \frac{d}{d r}-\mu_{i}^{2} r^{2}\right) \sigma_{i, 0}(r)&=0 &\text { in } \Omega_{i},\\
\left(r^{2} \frac{d^{2}}{d r^{2}}+r \frac{d}{d r}-\left(\mu_{i}^{2} r^{2}+l^{2}\right)\right) \sigma_{i, 1}(r)&=0 &\text { in } \Omega_{i}.
\end{align}Therefore, it is sufficient to consider the expression
\begin{align}
\sigma&=A_{1} I_{0}\left(\mu_{1} r\right)+\delta e^{i l \theta} B_{1} I_{l}\left(\mu_{1} r\right) &\text { in } \Omega_{1},\\
\sigma&=A_{2} I_{0}\left(\mu_{2} r\right)+A_{3} K_{0}\left(\mu_{2} r\right)+\delta e^{i l \theta}\left(B_{2} I_{l}\left(\mu_{2} r\right)+B_{3} K_{l}\left(\mu_{2} r\right)\right) &\text { in } \Omega_{2}.
\end{align}Apply $\eqref{noslipperybc},\eqref{fluxbc},\eqref{farfieldbc}$ on the interface $r=R+\delta e^{i l \theta}$ with $\delta\ll1$. Orders higher than $O(\delta)$ are all discarded in the following calculations.

At $O(1)$, the equations are the same as the radial solution.

The equations at $O(\delta)$ determine the coefficients $B_{1}, B_{2}, B_{3}:$
\begin{align}
B_{1} I_{l}\left(\mu_{1} R\right)=B_{2} I_{l}\left(\mu_{2} R\right)+B_{3} K_{l}\left(\mu_{2} R\right)+\left(A_{1} \mu_{1} I_{1}\left(\mu_{1} R\right)\right)\left(\frac{1}{D}-1\right),\\
0=B_{2} I_{l}\left(\mu_{2} R_{\infty}\right)+B_{3} K_{l}\left(\mu_{2} R_{\infty}\right),\\
\left(B_{1} \mu_{1}\left(I_{l-1}\left(\mu_{1} R\right)-\frac{l}{\mu_{1} R} I_{l}\left(\mu_{1} R\right)\right)+A_{1} \mu_{1}^{2}\left(I_{0}\left(\mu_{1} R\right)-\frac{1}{\mu_{1} R} I_{1}\left(\mu_{1} R\right)\right)\right)=\nonumber\\
D\left(B_{2} \mu_{2}\left(I_{l-1}\left(\mu_{2} R\right)-\frac{l}{\mu_{2} R} I_{l}\left(\mu_{2} R\right)\right)-B_{3} \mu_{2}\left(K_{l-1}\left(\mu_{2} R\right)+\frac{l}{\mu_{2} R} K_{l}\left(\mu_{2} R\right)\right)\right.\nonumber\\
+\left.A_{2} \mu_{2}^{2}\left(I_{0}\left(\mu_{2} R\right)-\frac{1}{\mu_{2} R} I_{1}\left(\mu_{2} R\right)\right)+A_{3} \mu_{2}^{2}\left(K_{0}\left(\mu_{2} R\right)+\frac{1}{\mu_{2} R} K_{1}\left(\mu_{2} R\right)\right)\right).
\end{align}
Solving for $B_{1}, B_{2}, B_{3}$, we have\\
\begin{equation}
\small
\begin{aligned}
B_1=\left({D}(I_{l}(R_{\infty} \mu_{2}) K_{l}(R \mu_{2})-I_{l}(R \mu_{2}) K_{l}(R_{\infty} \mu_{2})) \mu_{2}(A_{2}(R I_{0}(R \mu_{2}) \mu_{2}-I_{1}(R \mu_{2}))
+A_{3}(K_{1}(R \mu_{2})+R K_{0}(R \mu_{2}) \mu_{2}))\right.\\
+A_{1} \mu_{1}\left(R I_{0}(R \mu_{1})(I_{l}(R \mu_{2}) K_{l}(R_{\infty} \mu_{2})
-I_{l}(R_{\infty} \mu_{2}) K_{l}(R \mu_{2})) \mu_{1}
+I_{1}(R \mu_{1})\left((l({D}-1)-1) I_{l}(R \mu_{2}) K_{l}(R_{\infty} \mu_{2})\right.\right.\\
-\left.\left.\left.R({D}-1) I_{l-1}(R \mu_{2}) \mu_{2} K_{l}(R_{\infty} \mu_{2})
+I_{l}(R_{\infty} \mu_{2})\left(K_{l}(R \mu_{2})-R({D}-1) \mu_{2}\left(K_{l-1}(R \mu_{2})+\frac{l K_{l}(R \mu_{2})}{R \mu_{2}}\right)\right)\right)\right)\right)\\
/
\left(R {D} I_{l-1}(R \mu_{2}) I_{l}(R \mu_{1}) K_{l}(R_{\infty} \mu_{2}) \mu_{2}+R(I_{l}(R_{\infty} \mu_{2}) K_{l}(R \mu_{2})-I_{l}(R \mu_{2}) K_{l}(R_{\infty} \mu_{2})) \mu_{1}\left(I_{l-1}\left(R \mu_{1}\right)-\frac{l I_{l}(R \mu_{1})}{R_{\mu_{1}}}\right)\right.\\
+\left.{D} I_{l}\left(R \mu_{1}\right)\left(R I_{l}\left(R_{\infty} \mu_{2}\right) \mu_{2}\left(K_{l-1}\left(R \mu_{2}\right)+\frac{l K_{l}\left(R \mu_{2}\right)}{R \mu_{2}}\right)-l I_{l}\left(R \mu_{2}\right) K_{l}\left(R_{\infty} \mu_{2}\right)\right)\right)
,\nonumber
\end{aligned}
\end{equation}
\begin{equation}
\begin{aligned}
B_2=-\left(K_{l}(R_{\infty} \mu_{2})\left(I_{l}(R \mu_{1}) \mu_{2}(A_{2}(R I_{0}(R \mu_{2}) \mu_{2}-I_{1}(R \mu_{2}))+A_{3}(K_{1}(R \mu_{2})+R K_{0}\left(R \mu_{2}) \mu_{2})\right) {D}^{2}\right.\right.\\
\left.\left.+A_{1} \mu_{1}\left(I_{1}\left(R \mu_{1}\right)\left({D} I_{l}\left(R \mu_{1}\right)+R({D}-1) \mu_{1}\left(I_{l-1}\left(R \mu_{1}\right)-\frac{l I_{l}\left(R \mu_{1}\right)}{R \mu_{1}}\right)\right)-R {D} I_{0}\left(R \mu_{1}\right) I_{l}\left(R \mu_{1}\right) \mu_{1}\right)\right)\right)\\
/
\left({D}\left( R {D} I_{l-1}(R \mu_{2}) I_{l}(R \mu_{1}) K_{l}(R_{\infty} \mu_{2}) \mu_{2}+R(I_{l}(R_{\infty} \mu_{2}) K_{l}(R \mu_{2})-I_{l}(R \mu_{2}) K_{l}(R_{\infty} \mu_{2})) \mu_{1}
\left(I_{l-1}(R \mu_{1})-\frac{l I_{l}(R \mu_{1})}{R \mu_{1}}\right)\right.\right.\\
\left.\left.+{D} I_{l}(R \mu_{1})\left(R I_{l}(R_{\infty} \mu_{2}) \mu_{2}
\left(K_{l-1}(R \mu_{2})+\frac{l K_{l}(R \mu_{2})}{R \mu_{2}}\right)
-l I_{l}(R \mu_{2}) K_{l}(R_{\infty} \mu_{2})\right) \right)\right)
,\nonumber
\end{aligned}
\end{equation}
\begin{equation}
\begin{aligned}
B_3=
\left(I_{l}\left(R_{\infty} \mu_{2}\right)\left({D} I_{l}\left(R \mu_{1}\right)\left({D} \mu_{2}\left(A_{2}\left(R I_{0}\left(R \mu_{2}\right) \mu_{2}-I_{1}\left(R \mu_{2}\right)\right)+A_{3}\left(K_{1}\left(R \mu_{2}\right)+R K_{0}\left(R \mu_{2}\right) \mu_{2}\right)\right)\right.\right.\right.\\
\left.\left.\left.-R I_{0}\left(R \mu_{1}\right) A_{1} \mu_{1}^{2}\right)+I_{1}\left(R \mu_{1}\right) A_{1} \mu_{1}\left({D} I_{l}\left(R \mu_{1}\right)+R({D}-1) \mu_{1}\left(I_{l-1}\left(R \mu_{1}\right)-\frac{l I_{l}\left(R \mu_{1}\right)}{R \mu_{1}}\right)\right) \right)\right)\\
/\left({D}\left(R {D} I_{l-1}\left(R \mu_{2}\right) I_{l}\left(R \mu_{1}\right) K_{l}\left(R_{\infty} \mu_{2}\right) \mu_{2}+R\left(I_{l}\left(R_{\infty} \mu_{2}\right) K_{l}\left(R \mu_{2}\right)\right.\right.\right.\\
-\left.\left.I_{l}\left(R \mu_{2}\right) K_{l}\left(R_{\infty} \mu_{2}\right)\right) \mu_{1}\left(I_{l-1}\left(R \mu_{1}\right)-\frac{l I_{l}\left(R \mu_{1}\right)}{R \mu_{1}}\right)\right.\\
+\left.\left.{D} I_{l}\left(R \mu_{1}\right)\left(R I_{l}\left(R_{\infty} \mu_{2}\right) \mu_{2}\left(K_{l-1}\left(R \mu_{2}\right)+\frac{l K_{l}\left(R \mu_{2}\right)}{R \mu_{2}}\right)-l I_{l}\left(R \mu_{2}\right) K_{l}\left(R_{\infty} \mu_{2}\right)\right)\right)\right).\nonumber
\end{aligned}
\end{equation}
The nutrient $\sigma$ on $\Gamma$ is given by\\
\begin{eqnarray}{}
(\sigma)_{\Gamma}&=&\left(A_{1} I_{0}\left(\mu_{1} r\right)+B_{1} I_{l}\left(\mu_{1} r\right) \delta e^{i l \theta}\right)_{\Gamma}\nonumber\\
&=&A_{1} I_{0}\left(\mu_{1} R\right)+\left(\mu_{1} A_{1} I_{1}\left(\mu_{1} R\right)+B_{1} I_{l}\left(\mu_{1} R\right)\right) \delta e^{i l \theta}.\label{nutatbdry}
\end{eqnarray}\\
The normal derivative of $\sigma$ on $\Gamma$ is given by\\
\begin{eqnarray}{} \left(\frac{\partial \sigma}{\partial \mathbf{n}}\right)_{\Gamma} &=&\left(\frac{\partial \sigma}{\partial r}\right)_{\Gamma} \nonumber\\&=&\left(\left(A_{1} I_{0}\left(\mu_{1} r\right)+B_{1} I_{l}\left(\mu_{1} r\right) \delta e^{i l \theta}\right)_{r}\right)\nonumber \\&=&\left(A_{1} \mu_{1} I_{1}\left(\mu_{1} r\right)+\left(B_{1} \mu_{1}\left(I_{l-1}\left(\mu_{1} r\right)-\frac{l}{\mu_{1} r} I_{l}\left(\mu_{1} r\right)\right)\right) \delta e^{i l \theta}\right)_{\Gamma}\nonumber \\
&=&\mu_{1} A_{1} I_{1}\left(\mu_{1} R\right)+\left(A_{1}\left(\mu_{1}^{2} I_{0}\left(\mu_{1} R\right)-\mu_{1} \frac{I_{1}\left(\mu_{1} R\right)}{R}\right)+B_{1}\left(\mu_{1} I_{l-1}\left(\mu_{1} R\right)-l \frac{I_{l}\left(\mu_{1} R\right)}{R}\right)\right) \delta e^{i l \theta},\nonumber\\
\label{nutfluxatbdry}
\end{eqnarray}\\
where we use the identity in Eq. $\eqref{besseliid}$,
\begin{equation}
I_{2}\left(\mu_{1} R\right)=2 I_{1}\left(\mu_{1} R\right)^{\prime}-I_{0}\left(\mu_{1} R\right)=2\left(I_{0}\left(\mu_{1} R\right)-\frac{1}{\mu_{1} R} I_{1}\left(\mu_{1} R\right)\right)-I_{0}\left(\mu_{1} R\right).\nonumber
\end{equation}Similarly we seek a solution of Laplace equation, which reduces to Eq. $\eqref{pressurebc}$ on the perturbed circle given by Eq. $\eqref{perturbedcircle}$.
It is sufficient to consider the expression
\begin{equation}
p=A+\delta e^{i l \theta} B r^{l}.
\end{equation}For the perturbed circle defined by Eq. $\eqref{perturbedcircle}, \kappa$ is given by
\begin{equation}
\kappa=\frac{1}{R}\left(1+\frac{l^{2}-1}{R} \delta e^{i l \theta}\right).
\end{equation}
On the interface from Eqs. $\eqref{pressurebc},\eqref{perturbedcircle}$, and $\eqref{nutatbdry}$ we obtain
\begin{eqnarray}{}
(p)_{\Gamma}
&=&A+B R^{l} \delta e^{i l \theta} \nonumber\\
&=&\mathcal{G}^{-1} \frac{1}{R}+\left(\mathcal{P}-\chi_{\sigma}\right)\left(A_{1} I_{0}\left(\mu_{1} R\right)\right)-\frac{\mathcal{A}}{2 d} R^{2}\nonumber\\
&&+\left(\mathcal{G}^{-1} \frac{l^{2}-1}{R^{2}}-\frac{\mathcal{A}}{d} R+\left(\mathcal{P}-\chi_{\sigma}\right)\left(A_{1} \mu_{1} I_{1}\left(\mu_{1} R\right)+B_{1} I_{l}\left(\mu_{1} R\right)\right)\right) \delta e^{i l \theta}.
\end{eqnarray}
It is straightforward to derive that
\begin{align}{}
A&=\mathcal{G}^{-1} \frac{1}{R}-\frac{\mathcal{A}}{2 d} R^{2}+\left(\mathcal{P}-\chi_{\sigma}\right) A_{1} I_{0}\left(\mu_{1} R\right),\\
B&=\left(\mathcal{G}^{-1} \frac{l^{2}-1}{R^{2}}-\frac{\mathcal{A}}{d} R+\left(\mathcal{P}-\chi_{\sigma}\right)\left(\mu_{1} A_{1} I_{1}\left(\mu_{1} R\right)+B_{1} I_{l}\left(\mu_{1} R\right)\right)\right) \frac{1}{R^{l}}.
\end{align}
Thus
\begin{eqnarray}{}
p&=&\mathcal{G}^{-1} \frac{1}{R}-\frac{\mathcal{A}}{2 d} R^{2}+\left(\mathcal{P}-\chi_{\sigma}\right) A_{1} I_{0}\left(\mu_{1} R\right)\nonumber\\
&&+\left(\mathcal{G}^{-1} \frac{l^{2}-1}{R^{2}}-\frac{\mathcal{A}}{d} R+\left(\mathcal{P}-\chi_{\sigma}\right)\left(\mu_{1} A_{1} I_{1}\left(\mu_{1} R\right)+B_{1} I_{l}\left(\mu_{1} R\right)\right)\right) \frac{r^{l}}{R^{l}} \delta e^{i l \theta}.
\end{eqnarray}
The normal derivative of $p$ is given by
\begin{eqnarray}{} \left(\frac{\partial p}{\partial \mathbf{n}}\right)_{\Gamma} &=&\left(\frac{\partial p}{\partial r}\right)_{\Gamma}\nonumber\\
&=&l\left(\mathcal{G}^{-1} \frac{l^{2}-1}{R^{3}}-\frac{\mathcal{A}}{d}+\left(\mathcal{P}-\chi_{\sigma}\right)\left(\mu_{1} A_{1} \frac{I_{1}\left(\mu_{1} R\right)}{R}+B_{1} \frac{I_{l}\left(\mu_{1} R\right)}{R}\right)\right) \delta e^{i l \theta}\label{pressureatbdry}.
\end{eqnarray}Note that
\begin{equation}
n \cdot(\mathbf{x})_{\Gamma}=r+O\left(\delta^{2}\right)=R+\delta e^{i l \theta}+O\left(\delta^{2}\right).\label{normalxatbdry}
\end{equation}
Combining $\eqref{nutfluxatbdry},\eqref{pressureatbdry},\eqref{normalxatbdry}$ and $\eqref{normalvel}$, we obtain
\begin{eqnarray}{}
V&= &\frac{dR}{dt}+\frac{d\delta}{dt} e^{i l \theta} \nonumber\\
&=&-l\left(\mathcal{G}^{-1} \frac{l^{2}-1}{R^{3}}-\frac{\mathcal{A}}{d}+\left(\mathcal{P}-\chi_{\sigma}\right)\left(\mu_{1} A_{1} \frac{I_{1}\left(\mu_{1} R\right)}{R}+B_{1} \frac{I_l\left(\mu_{1} R\right)}{R}\right)\right) \delta e^{i l \theta}+\mathcal{P} \frac{\partial \sigma}{\partial r}-\frac{\mathcal{A}}{d}\left(R+\delta e^{i l \theta}\right) \nonumber\\
&=&-\frac{\mathcal{A} R}{d}+\mathcal{P} \frac{\partial \sigma}{\partial r}+\left(-\mathcal{G}^{-1} \frac{l\left(l^{2}-1\right)}{R^{3}}+\frac{l-1}{d} \mathcal{A}-l\left(\mathcal{P}-\chi_{\sigma}\right)\left(\mu_{1} A_{1} \frac{I_{1}\left(\mu_{1} R\right)}{R}+B_{1} \frac{I_{l}\left(\mu_{1} R\right)}{R}\right)\right) \delta e^{i l \theta}\nonumber\\
&=&- \frac{\mathcal{A} R}{d}+\mathcal{P}\left(\mu_{1} A_{1} I_{1}\left(\mu_{1} R\right)
+\left(A_{1}\left(\mu_{1}^{2} I_{0}\left(\mu_{1} R\right)-\mu_{1} \frac{I_{1}\left(\mu_{1} R\right)}{R}\right)+B_{1}\left(\mu_{1} I_{l-1}\left(\mu_{1} R\right)-l \frac{I_{1}\left(\mu_{1} R\right)}{R}\right)\right) \delta e^{i l \theta}\right)\nonumber\\
&&+\left(-\mathcal{G}^{-1} \frac{l\left(l^{2}-1\right)}{R^{3}}+\frac{l-1}{d} \mathcal{A}-l\left(\mathcal{P}-\chi_{\sigma}\right)\left(\mu_{1} A_{1} \frac{I_{1}\left(\mu_{1} R\right)}{R}+B_{1} \frac{I_l\left(\mu_{1} R\right)}{R}\right)\right) \delta e^{i l \theta} \nonumber\\
&=&C \mathcal{P}-\frac{\mathcal{A} R}{d}+\nonumber\\
&&\left(-\mathcal{G}^{-1} \frac{l\left(l^{2}-1\right)}{R^{3}}+\frac{l-1}{d} \mathcal{A}+ \mathcal{P}\left(\mu_{1}^{2} A_{1} I_{0}\left(\mu_{1} R\right)+B_{1}\left(\mu_{1} I_{l-1}\left(\mu_{1} R\right)-\frac{l}{R} I_{l}\left(\mu_{1} R\right)\right)-\frac{{C}}{R}\right)\right.\nonumber\\
&&-\left.\left(\mathcal{P}-\chi_{\sigma}\right)\left(B_{1} \frac{l}{R} I_{l}\left(\mu_{1} R\right)+\frac{l}{R} {C}\right)\right) \delta e^{i l \theta} ,
\end{eqnarray}
where $C=\mu_{1} A_{1} I_{1}\left(\mu_{1} R\right)$.\\\\
\noindent
Equating coefficients of like harmonics, we obtain
\begin{equation}
\frac{dR}{dt}=C \mathcal{P}-\frac{\mathcal{A} R}{d},
\end{equation}
\begin{equation}
R^{-1}\frac{dR}{dt}=\frac{C \mathcal{P}}{R}-\frac{\mathcal{A}}{d},
\end{equation}
and
\begin{eqnarray}{}
\delta^{-1}\frac{d\delta}{dt}&=&-\mathcal{G}^{-1} \frac{l\left(l^{2}-1\right)}{R^{3}}+\frac{l-1}{d} \mathcal{A}\nonumber\\
&&+ \mathcal{P}\left(\mu_{1}^{2} A_{1} I_{0}\left(\mu_{1} R\right)+B_{1}\left(\mu_{1} I_{l-1}\left(\mu_{1} R\right)-\frac{l}{R} I_{l}\left(\mu_{1} R\right)\right)-\frac{{C}}{R}\right)\nonumber\\
&&-\left(\mathcal{P}-\chi_{\sigma}\right)\left(\mu_{1} B_{1} \frac{l}{R} I_{l}\left(\mu_{1} R\right)+\frac{l}{R} {C}\right).
\end{eqnarray}
The equation of shape perturbation is given by
\begin{eqnarray}{}
\left(\frac{\delta}{R}\right)^{-1}
\frac{d}{dt}{\left(\frac{\delta}{R}\right)}
&=&\delta^{-1}\frac{d\delta}{dt}-R^{-1}\frac{dR}{dt}\nonumber\\
&=&-\mathcal{G}^{-1} \frac{l\left(l^{2}-1\right)}{R^{3}}+\frac{l}{d} \mathcal{A}\nonumber\\
&&+\mathcal{P}\left(\mu_{1}^{2} A_{1} I_{0}\left(\mu_{1} R\right)+B_{1}\left(\mu_{1} I_{l-1}\left(\mu_{1} R\right)-\frac{l}{R} I_{l}\left(\mu_{1} R\right)\right)-\frac{2}{R} {C}\right)\nonumber\\
&&-\left(\mathcal{P}-\chi_{\sigma}\right)\left(B_{1} \frac{l}{R} I_{l}\left(\mu_{1} R\right)+\frac{l}{R} C\right).
\end{eqnarray}
The critical apoptosis parameter $\mathcal{A}_{c}$ is the function of $R$ such that $\frac{d}{dt}{\left(\frac{{\delta}}{R}\right)}=0$ and is given by
\begin{eqnarray}{}
\mathcal{A}_{c}= \mathcal{G}^{-1} \frac{d\left(l^{2}-1\right)}{R^{3}}
&&-\mathcal{P} \frac{d}{l}\left(\mu_{1}^{2} A_{1} I_{0}\left(\mu_{1} R\right)+B_{1}\left(\mu_{1} I_{l-1}\left(\mu_{1} R\right)-\frac{l}{R} I_{l}\left(\mu_{1} R\right)\right)-\frac{2}{R} {C}\right)\nonumber\\ &&+\left(\mathcal{P}-\chi_{\sigma}\right) \frac{d}{l}\left(B_{1} \frac{l}{R} I_{l}\left(\mu_{1} R\right)+\frac{l}{R} {C}\right).
\end{eqnarray}

\section*{Appendix B: The Evaluation of the Boundary Integrals}
\label{Evaluation of BIM}
With the integral formulation above, we assume interface curves $\Gamma$ and $\Gamma_{\infty}$ are analytic and given by $\big\{\mathbf{x}(\alpha,t)=(x(\alpha,t),y(\alpha,t): 0\leq \alpha \leq 2 \pi \big\}$, where $\mathbf{x}$ is $2 \pi$-periodic in the parametrization $\alpha$. The unit tangent and normal(outward) vectors can be calculated as $\mathbf{s}=(x_\alpha,y_\alpha)/s_\alpha$, $\mathbf{n}=(y_\alpha,-x_\alpha)/s_\alpha$, where the local variation of the arclength $s_\alpha=\sqrt{x_\alpha^2+y_\alpha^2}$. Subscripts refer to partial differentiation.
We track the interfaces $\Gamma$ and $\Gamma_{\infty}$ by introducing N marker points to discretize the planar curves, parametrized by $\alpha_j=jh$, $h=\frac{2\pi}{N}$, $N$ is a power of $2$. Here we focus on the numerical evaluation of integrals following \cite{jou1997,li2011boundary,minjhe2019}. A rigorous convergence and error analysis of the boundary integral method for a simplified tumor problem can be found in \cite{Wenrui2018}.

\paragraph{Computation of the single-layer potential type integral.\\\\}
In Eqs. \eqref{eq40}, \eqref{eq41}, \eqref{eq42} and \eqref{eq49}, the single-layer potential type integrals contain the Green functions with a logarithmic singularity at $r=0$. They can be rewritten in the following form under the parametrization $\alpha$
\begin{equation}\label{eq50}
    \int_{\Gamma}
    \Phi(\alpha,\alpha')\phi(\alpha')s_{\alpha}(\alpha')d\alpha',
\end{equation}where $\Phi$ are the Green functions $G$ or $G_i$ , $\Gamma$ may be either $\Gamma$ or $\Gamma_{\infty}$ and $\phi$ may be $\eta$ ,$\frac{\partial\sigma_1}{\partial\mathbf{n}}$ or $\frac{\partial\sigma_2}{\partial\mathbf{n}_\infty}$. We may decompose the Green functions as below
\begin{equation}\label{eq51}
    G(\alpha,\alpha')=
    -\frac{1}{2 \pi} \ln r
    =-\frac{1}{2\pi}\left(
    \ln{2 \left| \sin{\frac{\alpha-\alpha'}{2}}\right|}
    +\left[\ln r-\ln{2 \left| \sin{\frac{\alpha-\alpha'}{2}}\right|}\right]
    \right),
\end{equation}

\begin{equation}\label{eq52}
    G_{i}(\alpha,\alpha')=
    \frac{1}{2 \pi} K_{0}(\mu_{i} r)=
    -\frac{1}{2 \pi}
    \left(I_{0}(\mu_{i}r)\ln{2 \left| \sin{\frac{\alpha-\alpha'}{2}}\right|}
    +\left[-K_{0}(\mu_{i}r)
    -I_{0}(\mu_{i}r)\ln{2 \left| \sin{\frac{\alpha-\alpha'}{2}}\right|}
    \right]\right),
\end{equation}where $I_{0}$ is a modified Bessel function of the first kind, $r=|\mathbf{x}(\alpha)-\mathbf{x}'(\alpha ')|$. The square brackets on the right-hand side of Eqs.\eqref{eq51}, \eqref{eq52} have removable singularity at $\alpha=\alpha'$, since $r=
s_{\alpha}\left|\alpha-\alpha'\right|
\sqrt{1+\mathscr{O}(\alpha-\alpha')}
=s_{\alpha}\left|\alpha-\alpha'\right|
(1+\mathscr{O}(\alpha-\alpha'))$ for $\alpha \approx \alpha'$, where $\mathscr{O(\alpha-\alpha')}$ denotes a smooth function that vanishes as $\alpha\rightarrow\alpha'$, and since $K_0$ has the expansion
\begin{equation}
K_{0}(z)=-\left(\log \frac{z}{2}+C\right) I_{0}(z)+\Sigma_{n=1}^{\infty} \frac{\psi(n)}{(n !)^{2}}\left(\frac{z}{2}\right)^{2 n}.
\end{equation}
Thus, for an analytic and $2\pi$-periodic function $f(\alpha,\alpha')$, a standard trapezoidal rule or alternating point rule can be used to evaluate the integral
\begin{equation}\label{eq53}
    \int_{0}^{2\pi}
f(\alpha,\alpha')
\ln{\frac{r}{2 \left| \sin{\frac{\alpha-\alpha'}{2}}\right|}}
d\alpha'.
\end{equation}The remaining terms on the right-hand side of Eqs.\eqref{eq51}, \eqref{eq52} have logarithmic singularity and can be evaluated through the following spectrally accurate quadrature \cite{kress1995numerical}
\begin{equation}\label{eq54}
    \int_{0}^{2\pi}f(\alpha_i,\alpha')
\ln{2 \left| \sin{\frac{\alpha_i-\alpha'}{2}}\right|}
d\alpha'\approx
\Sigma_{j=0}^{2m-1}q_{\left|j-i\right|}f(\alpha_i,\alpha_j),
\end{equation}where $m=\frac{N}{2}$, $\alpha_i=\frac{\pi i}{m}$ for $i=0,1,...,2m-1$, and weight coefficients
\begin{equation}\label{eq55}
    q_j=-\frac{\pi}{m}\Sigma_{k=1}^{m-1}\frac{1}{k}\cos{\frac{kj\pi}{m}}-\frac{(-1)^j\pi}{2m^2} , \text{for } j=0,1,...,2m-1.
\end{equation}The derivative $\frac{d}{ds}$ in Eq. \eqref{eq49} is approximated using fast Fourier transform spectral derivatives thus maintaining spectral accuracy.
\paragraph{Computation of the double-layer potential-type integral\\\\}
In Eqs. \eqref{eq40}, \eqref{eq41}, \eqref{eq42} and \eqref{eq47}, the double-layer potential type integrals contain the Green functions with singularity at $r=0$ (logarithmic for Eqs.\eqref{eq40}, \eqref{eq41}, \eqref{eq42}). They can be rewritten as in the following form under the parametrization $\alpha$
\begin{equation}\label{eq56}
    \int_{\Gamma}
    \frac{\partial \Phi(\alpha,\alpha')}{\partial \mathbf{n}(\alpha')}\phi(\alpha')s_{\alpha}(\alpha')d\alpha',
\end{equation}where $\Phi$ are the Green functions $G$ or $G_i$, $\Gamma$ may be either $\Gamma$ or $\Gamma_{\infty}$ and $\phi$ may be $\eta$, $\sigma_1$ or $\sigma_2=1$. Further,
\begin{equation}\label{eq57}
    \frac{\partial G(\alpha,\alpha')}{\partial \mathbf{n}(\alpha')}s_{\alpha}(\alpha')=
    h(\alpha,\alpha')\frac{1}{r},
\end{equation}where the auxiliary function $h(\alpha,\alpha')=\frac{(\mathbf{x(\alpha)}-\mathbf{x(\alpha')})\cdot\mathbf{n(\alpha')}s_\alpha(\alpha')}{2\pi r}$ with $r=\left|\mathbf{x(\alpha)}-\mathbf{x(\alpha')}\right|$. Note that $h(\alpha,\alpha')\sim\mathscr{O}(\alpha-\alpha')$.
Since $\frac{\partial G}{\partial \mathbf{n}}$ has no logarithmic singularity, we may simply use the alternating point rule to evaluate it.
For $\frac{\partial G_i}{\partial \mathbf{n}}$, we decompose it as below
\begin{equation}\label{eq58}
    \frac{\partial G_{i}(\alpha,\alpha')}{\partial \mathbf{n}(\alpha')}s_{\alpha}(\alpha')=
    h(\alpha,\alpha') K_{1}(\mu_{i} r)=
    g_1(\alpha,\alpha')\ln{2 \left| \sin{\frac{\alpha-\alpha'}{2}}\right|}
    +g_2(\alpha,\alpha'),
\end{equation}where $g_1(\alpha,\alpha')$ and $g_2(\alpha,\alpha')$ are analytic and $2\pi$-periodic functions with
\begin{equation}\label{eq59}
    g_1(\alpha,\alpha')=h(\alpha,\alpha')I_1(\mu_i r),
\end{equation}
\begin{equation}\label{eq60}
    g_2(\alpha,\alpha')=h(\alpha,\alpha')
    \left[K_{1}(\mu_{i}r)
    -I_{1}(\mu_{i}r)\ln{2 \left| \sin{\frac{\alpha-\alpha'}{2}}\right|}
    \right],
\end{equation}
where we have used the fact
\begin{equation}
\frac{d}{dr}K_0(r)=-K_1(r).
\end{equation}
Since $K_1$ has the expansion
\begin{equation}
K_{1}(z)=\frac{1}{z}+\left(\log \frac{z}{2}+C\right) I_{1}(z)-\frac{1}{2} \sum_{n=0}^{\infty} \frac{\psi(n+1)+\psi(n)}{n !(n+1) !}\left(\frac{z}{2}\right)^{2 n+1},
\end{equation}
the  square bracket on the right-hand side of Eq. $\eqref{eq60}$ also has removable singularity at $\alpha=\alpha'$ thus the integral involving $g_2(\alpha,\alpha')$ can be evaluated by a standard trapezoidal rule or alternating point rule. Note that
\begin{equation}\label{eq61}
    g_2(\alpha,\alpha)=\frac{h(\alpha,\alpha)}{\mu_i r}=\frac{1}{4\pi\mu_i}
    \frac{x_\alpha y_{\alpha\alpha}-x_{\alpha\alpha}y_\alpha}{x_\alpha^2+y_\alpha^2}.
\end{equation}The first term on the right-hand side of Eq.$\eqref{eq58}$ is still singular and evaluated through the quadrature given in Eqs. \eqref{eq54} and \eqref{eq55}.\\
To summarize, using $\text{Nystr\"om}$ discretization with the Kress quadrature rule described above, we reduce the boundary integral Eqs. \eqref{eq40}, \eqref{eq41}, \eqref{eq42} and  \eqref{eq47} to two dense linear systems with the unknowns as the discretization of $\eta$ , $\sigma_1$, $\frac{\partial\sigma_1}{\partial \mathbf{n'}}$ on $\Gamma$ and $\frac{\partial\sigma_2}{\partial\mathbf{n'_\infty}}$ on $\Gamma_\infty$, which can be solved using an iterative solver, e.g., GMRES \cite{saad1986gmres}.

\section*{Appendix C: The Evolution of the Interface}
\label{Evolution of the interface}
As indicated by \cite{hou1994removing}, the curvature-driven motion introduces high-order derivatives, both non-local and nonlinear, into the dynamics through the Laplace-Young condition at the interface. Explicit time integration methods thus suffer from severe stability constraints and implicit methods are difficult to apply since the stiffness enters nonlinearly. Hou et al. resolved these difficulties by adopting the $\theta-L$ formulation and the small-scale decomposition (SSD), which we apply here.
\paragraph{$\theta-L$ formulation.\\\\}
This formulation helps to circumvent the problem of point clustering. Consider a point $\mathbf{x}(\alpha,t)=(x(\alpha,t),y(\alpha,t))\in \Gamma(t)$. Denote the unit tangent and normal (outward) vectors as $\hat{\mathbf{s}}=(x_\alpha,y_\alpha)/s_\alpha$ and $\hat{\mathbf{n}}=(y_\alpha,-x_\alpha)/s_\alpha$, the normal velocity and tangent velocity by $V(\alpha,t)=u\cdot\hat{\mathbf{n}}$ and $T(\alpha,t)=u\cdot\hat{\mathbf{s}}$, respectively, where $u=\mathbf{x}_t=V \hat{\mathbf{n}}+T \hat{\mathbf{s}}$ gives the motion of $\Gamma(t)$. The tangent angle that the planar curve $\Gamma(t)$ forms with the horizontal axis at $\mathbf{x}$, called $\theta$, satisfies $\theta=\tan^{-1}{\frac{y_\alpha}{x_\alpha}}$. The length of one period of the curve is $L(t)=\int_0^{2\pi}s_\alpha d\alpha$, where $s_\alpha$, the derivative of the arclength, satisfies $s_\alpha^2=x_\alpha^2+y_\alpha^2$. Differentiating these two equations in time, we obtain the following evolution equations:
\begin{equation}\label{eq62}
    \theta_t=\kappa T - V_s=\frac{1}{s_\alpha}(\theta_\alpha T- V_\alpha),
\end{equation}
\begin{equation}\label{eq63}
    s_{\alpha t}=(T_s+\kappa V)s_\alpha=T_\alpha+\theta_\alpha V.
\end{equation}
Instead of using the $(x,y)$ coordinates, $(L,\theta)$ becomes the dynamical variables. The unit tangent and normal vectors become $\Hat{\mathbf{s}}=(\cos{\theta},\sin{\theta})$, $\Hat{\mathbf{n}}=(\sin{\theta},-\cos{\theta})$.

 The normal velocity $V$ is calculated using Eq. \eqref{eq32}. The tangent velocity $T$ is chosen (independent of the morphology of the interface) such that the marker points are equally spaced in arclength to prevent point clustering:
\begin{equation}\label{eq64}
    T(\alpha,t)= \frac{\alpha}{2\pi}\int_0^{2\pi}\theta_{\alpha'}V' d\alpha'-\int_0^\alpha \theta_{\alpha'}V' d\alpha'.
\end{equation}
It follows that $s_\alpha$ is independent of $\alpha$ thus is everywhere equal to its mean:
\begin{equation}\label{eq65}
    s_\alpha=\frac{1}{2\pi}\int_0^{2\pi}s_\alpha(\alpha,t)d\alpha=\frac{L(t)}{2\pi}.
\end{equation}
The procedure for obtaining the initial equal arclength parametrization is presented in "Appendix B" of \cite{baker1990connection}. The
idea is to solve the nonlinear equation
\begin{equation}\label{eq66}
    \int_0^{\alpha_j} s_{\beta}d\beta=\frac{j}{N}L
\end{equation}for $\alpha_j$ using Newton's method and evaluate the equal arclength marker points $\mathbf{x}(\alpha_j)$by interpolation in Fourier space.
We may recover the interface by simply integrating:
\begin{equation}\label{eq67}
    \mathbf{x}_\alpha=\mathbf{x}_s s_\alpha=\frac{L(t)}{2\pi}(\cos{\theta(\alpha,t)},\sin{\theta(\alpha,t)}).
\end{equation}
\paragraph{Small scale decomposition (SSD).\\\\}
The idea of the small scale decomposition (SSD) is to extract the dominant part of the equations at small spatial scales \cite{hou1994removing}. To remove the stiffness, we use SSD in our problem and develop an explicit, non-stiff time integration algorithm.  In Eqs. \eqref{eq40}, \eqref{eq41}, \eqref{eq42}, \eqref{eq47} and \eqref{eq49}, based on  the analysis of the single-layer- and double-layer- type terms, the only singularity in the integrands comes from the logarithmic kernel. Following \cite{hou1994removing} and noticing the curvature term in Eq. \eqref{eq47}, one can show that at small spatial scales,
\begin{equation}\label{eq68}
V(\alpha,t) \sim \frac{1}{s_\alpha^2} \mathcal{H}[\theta_{\alpha\alpha}],
\end{equation}
where $\mathcal{H}(\xi)=\frac{1}{2\pi}\int_0^{2\pi}\xi'\cot{\frac{\alpha-\alpha'}{2}}d\alpha'$ is the Hilbert transform for a $2\pi$-periodic function $\xi$.\\ We rewrite Eq. \eqref{eq62},
\begin{equation}\label{eq69}
\theta_t=\frac{1}{s_\alpha^3} \mathcal{H}[\theta_{\alpha\alpha\alpha}]+N(\alpha,t),
\end{equation}
where the Hilbert transform term is the dominating high-order term at small spatial scales, and $ \displaystyle N= (\kappa T-V_s)-\frac{1}{s_\alpha^3} \mathcal{H}[\theta_{\alpha\alpha\alpha}]$ contains other lower-order terms in the evolution. This demonstrates that an explicit time-stepping method has the high-order constraint $\displaystyle \Delta t \le \left ( \frac{h}{s_\alpha} \right)^3$ where $\Delta t$ and $h$ are the time-step and spatial grid size, respectively. This has been demonstrated numerically in the seminal work \cite{hou1994removing} for a Hele-Shaw problem.  For the tumor growth problem, the semi-implicit time-stepping scheme (see Eq. \eqref{eq69}) requires $\Delta t = O(h)$ instead of explicit schemes which would require $\Delta t = O(h^3)$.  In Sect. \ref{sec:4.1}, we show a numerical example using $N=1024$ to simulate a 2-fold tumor. In this simulation, we could use $\Delta t$ as large as $\Delta t=1.0 \times 10^{-2}$ for stability instead of $\Delta t<10^{-6}$ for an explicit scheme with the equal-arclength parametrization. For the purpose of numerical accuracy, we used a smaller time step in our simulation.

\section*{Appendix D: Semi-implicit Time-Stepping Scheme}
\label{Time stepping}
Taking the Fourier transform of Eq. \eqref{eq69}, we get
\begin{equation}\label{eq70}
{\hat \theta}_t=-\frac{|k|^3}{s_\alpha^3} {\hat \theta}(k,t) +{\hat N}(k,t).
\end{equation}
We solve Eq. \eqref{eq70} using the second order accurate linear propagator method in the Adams-Bashforth form \cite{hou1994removing} in Fourier space and apply the inverse Fourier transform to recover $\theta$.  Specifically,  we discretize Eq. \eqref{eq70} as
\begin{equation}\label{eq71}
\quad\quad\;{\hat \theta}^{n+1}(k)=e_k(t_n,t_{n+1}){\hat \theta}^{n}(k)+\frac{\Delta t}{2}(3e_k(t_n,t_{n+1}){\hat N}^{n}(k)-e_k(t_{n-1},t_{n+1}){\hat N}^{n-1}(k),
\end{equation}
where  the superscript $n$ denotes the numerical solutions at $t=t_n$ and the integrating factor
\begin{equation}\label{eq72}
e_k(t_1,t_2)=\exp\left (-{|k|^3}\int_{t_1}^{t_2}\frac{dt}{s_\alpha^3(t)}\right ).
\end{equation}
Note that by setting the integrating factors in Eq. \eqref{eq71} to $1$, we recover the Adams-Bashforth explicit time-stepping method.
The integrating factors in Eq. \eqref{eq71} can be evaluated simply using the trapezoidal rule,
\begin{eqnarray}\label{eq73}
\int_{t_n}^{t_{n+1}}\frac{dt}{s_\alpha^3(t)} &\approx& \frac{\Delta t}{2} \left (\frac{1}{(s_\alpha^n)^3}+\frac{1}{(s_\alpha^{n+1})^3} \right ) \nonumber, \\
\int_{t_{n-1}}^{t_{n+1}}\frac{dt}{s_\alpha^3(t)} &\approx& {\Delta t} \left (\frac{1}{2(s_\alpha^{n-1})^3}+\frac{1}{(s_\alpha^{n})^3}+\frac{1}{2(s_\alpha^{n+1})^3} \right ).
\end{eqnarray}
To compute the arclength $s_\alpha$, equation \eqref{eq63} is discretized using the explicit second-order Adams-Bashforth method \cite{hou1994removing},
\begin{equation}\label{eq74}
s_\alpha^{n+1}=s_\alpha^n+\frac{\Delta t}{2}(3M^n-M^{n-1}),
\end{equation}
where $M$ is calculated using
\begin{equation}\label{eq75}
M=\frac{1}{2\pi}\int_0^{2\pi}V(\alpha,t)\theta_\alpha d\alpha.
\end{equation}

 Note that the second order linear propagator and Adams-Bashforth methods are multi-step method and require two previous time steps. The first time step is realized using an explicit Euler method for $s_\alpha^1$ and a first-order linear propagator of a similar form for $\hat{\theta}^1$.

 To reconstruct the tumor-host interface $(x(\alpha,t_{n+1}),y(\alpha,t_{n+1}))$ from the updated $\theta^{n+1}(\alpha)$ and $s_\alpha^{n+1}$, we first update a reference point $(x(0,t_{n+1}),y(0,t_{n+1})$ using a second-order explicit Adams-Bashforth method to discretize the equation of motion $\mathbf{x}_t=V\hat{\mathbf{n}}$ with the tangential part dropped since it does not change the morphology:
 \begin{equation}\label{eq76}
     (x(0,t_{n+1}),y(0,t_{n+1}))=(x(0,t_{n}),y(0,t_{n}))+\frac{\Delta t}{2} \left( 3V(0,t_n)\hat{\mathbf{n}}(0,t_{n})-V(0,t_{n-1})\hat{\mathbf{n}}(0,t_{n-1}) \right).
 \end{equation}
 Once we update the reference point, we obtain the configuration of the interface from the $\theta^{n+1}(\alpha)$ and $s_\alpha^{n+1}$ by integrating Eq. \eqref{eq67} following \cite{hou1994removing}:
 \begin{eqnarray}\label{eq77}
     x(\alpha,t_{n+1})&=&x(0,t_{n+1})+s_\alpha^{n+1}\left( \int_0^{\alpha} \cos(\theta^{n+1}(\alpha'))d\alpha'
     -\frac{\alpha}{2\pi}\int_0^{2\pi}\cos(\theta^{n+1}(\alpha'))d\alpha'\right),\nonumber\\
     y(\alpha,t_{n+1})&=&y(0,t_{n+1})+s_\alpha^{n+1}\left( \int_0^{\alpha} \sin(\theta^{n+1}(\alpha'))d\alpha'
     -\frac{\alpha}{2\pi}\int_0^{2\pi}\sin(\theta^{n+1}(\alpha'))d\alpha'\right),
 \end{eqnarray}where the indefinite integration is performed using the discrete Fourier transform.

 We use a 25th order Fourier filter to damp the highest nonphysical mode and suppress the  aliasing error \cite{hou1994removing}. We also use Krasny filtering \cite{krasny1986study} to prevent the accumulation of round-off errors during the computation.

We solve first the nutrient field $\sigma$ then the pressure field $p$. Next we compute the normal velocity $V$ and update the interface $\Gamma(t)$ and repeat this procedure.

\section*{Appendix E}
\label{D=1 and chi=10}

\begin{figure}
\includegraphics[width=\textwidth]{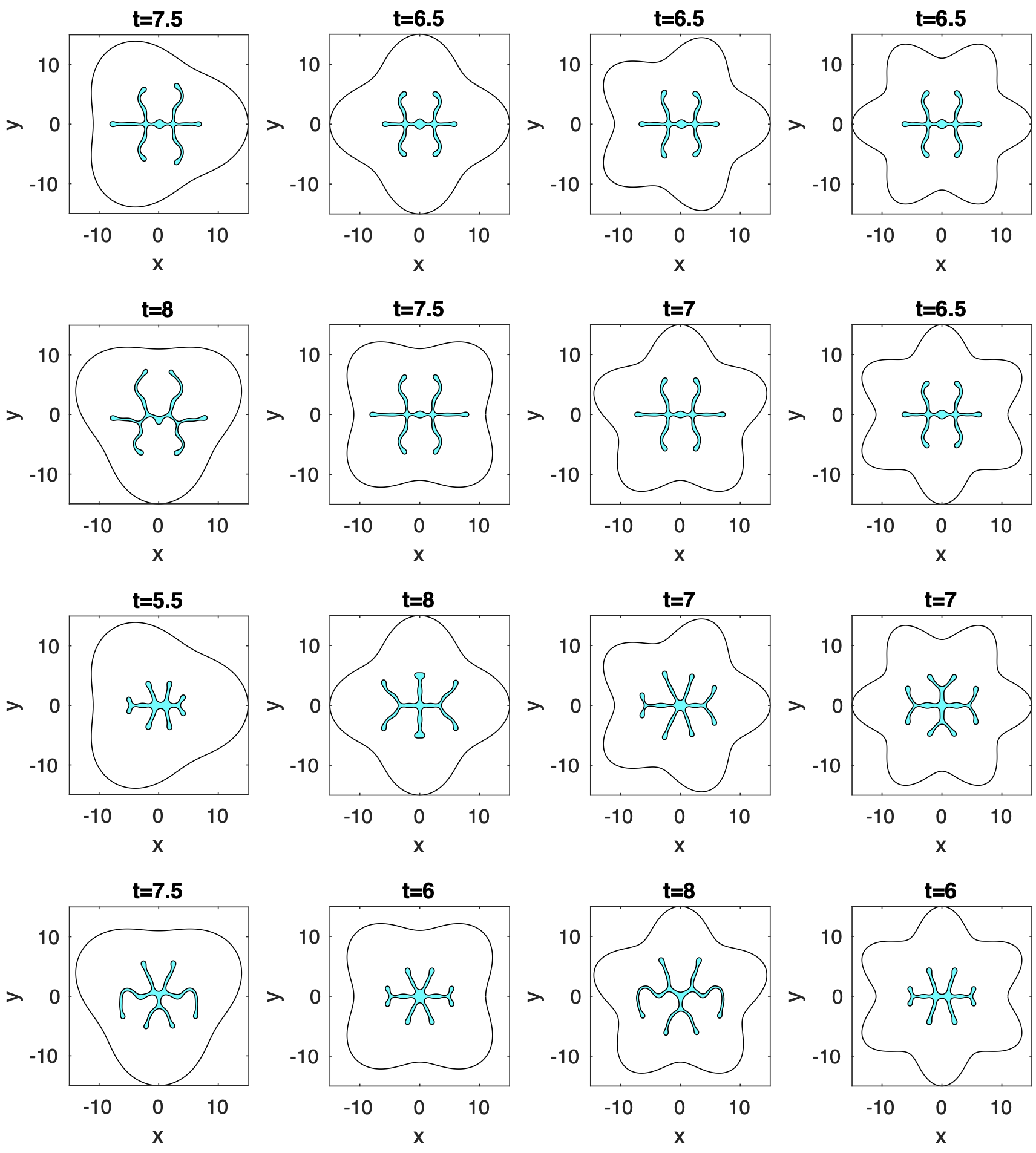}
\caption{Tumor morphologies with different symmetric far-field geometries in the nutrient-poor regime ($D=1$). The far-field boundaries are $R_\infty=13+2\cos(k\theta),k=3,4,5,6$ (Rows 1,3);
$R_\infty=13+2\cos(k\theta-\pi/2),k=3,5,~R_\infty=13+2\cos(k\theta-\pi),k=4,6$ (Rows 2,4). The initial tumor boundaries are $r=2.0+ 0.1\cos (2 \theta)$ (Rows 1,2); and $\frac{x^2}{2.1^2}+\frac{y^2}{1.9^2}=1$ (Rows 3,4). The remaining parameters are $\lambda=0.01$, $\chi_{\sigma}=10$, $\mathcal{P}=0.5$, $\mathcal{A}=0$,  $\mathcal{G}^{-1}=0.001$. Here, $N=512$, and $\Delta t=0.005$.
}
\label{fig:6-2}
\end{figure}

In Fig. \ref{fig:6-2} we present tumor morphologies at similar sizes under the same growth conditions but using  different symmetric far-field boundaries: $R_\infty=13+2\cos(k\theta),k=3,4,5,6$ (Row 1,3) and $R_\infty=13+2\cos(k\theta-\pi/2),k=3,5,R_\infty=13+2\cos(k\theta-\pi),k=4,6$ (Row 2,4). The parameters are the same as in Fig. \ref{fig:6} in the main text, where $\chi_\sigma=5$,  except that here in Fig. \ref{fig:6-2} we use $\chi_\sigma=10$. As in Fig. \ref{fig:6}, initial tumor boundary in rows 1 and 2 is the perturbed circle $R_\infty=2.0+ 0.1\cos (2 \theta)$ while in rows 3 and 4 the initial tumor boundary is the ellipse $\frac{x^2}{2.1^2}+\frac{y^2}{1.9^2}=1$. Hence, by comparing the evolution within each row, and between rows 1 and 2 and rows 3 and 4, we can see the effect of the far-field boundary shapes. By comparing rows 1 and  3 and rows 2 and 4, we can see the effect of the different initial shapes.

The results are similar to those obtained in Fig. \ref{fig:6} although here, because the chemotaxis coefficient $\chi_\sigma$ is increased, the tubular structures develop faster and are narrower than those observed in Fig. \ref{fig:6}. In particular, when the initial condition is the perturbed circle (rows 1 and 2), the far-field geometry has limited influence on tumor morphologies consistent with that observed in Figs. \ref{fig:4} and \ref{fig:5}. However, when the initial condition is an ellipse, which contains many modes, the morphologies are much more sensitive to the far-geometry because of the instability.

\bibliographystyle{spr-chicago}      
\bibliography{ComplexGeometryV2}   

\end{document}